\documentclass[12pt]{amsart}

\usepackage[top=1in,bottom=1in, left=1 in, right= 1 in, includehead,includefoot]{geometry}
\usepackage[T1]{fontenc}

\usepackage{hyperref,color}
\hypersetup{
     colorlinks=true,
     linkcolor=blue,
     filecolor=black,
     citecolor = blue,      
     urlcolor=blue,
     }

\usepackage{amsthm}
\usepackage{amsmath,amssymb}

\usepackage{mathtools} 
\usepackage[numbers,sort&compress]{natbib}
\usepackage{enumerate}

\newtheorem{theorem}{Theorem}[section]

\newtheorem{lemma}[theorem]{Lemma}
\newtheorem{proposition}[theorem]{Proposition}

\theoremstyle{definition}

\theoremstyle{remark}
\newtheorem{remark}[theorem]{Remark}
\theoremstyle{remark}

\numberwithin{equation}{section}

\newcommand{\ind}{{\bf 1}}

\def\inddd#1{{\ind}_{\left\{#1\right\}}} 
\newcommand{\proba}{\mathbb P}
\newcommand{\esp}{{\mathbb E}}

\newcommand{\inv}{^{-1}}

\newcommand{\eqnh}{\begin{eqnarray*}}
	\newcommand{\eqne}{\end{eqnarray*}}
\newcommand{\eqnhn}{\begin{eqnarray}}
	\newcommand{\eqnen}{\end{eqnarray}}
\newcommand{\equh}{\begin{equation}}
	\newcommand{\eque}{\end{equation}}

\def\summ#1#2#3{\sum_{#1 = #2}^{#3}}
\def\prodd#1#2#3{\prod_{#1 = #2}^{#3}}

\newcommand{\eqd}{\stackrel{d}{=}}

\def\topp#1{^{(#1)}}

\def\nn#1{{\left\|#1\right\|}}

\def\abs#1{\left|#1\right|}

\def\ccbb#1{\left\{#1\right\}} 

\def\pp#1{\left(#1\right)}
\def\spp#1{(#1)}
\def\bb#1{\left[#1\right]}

\def\mmid{\;\middle\vert\;}

\def\qmand{\quad\mbox{ and }\quad}

\def\qmwith{\quad\mbox{ with }\quad}
\def\mfa{\mbox{ for all }}

\def\wt#1{\widetilde{#1}}
\def\wb#1{\overline{#1}}

\def\limn{\lim_{n\to\infty}}

\def\limsupn{\limsup_{n\to\infty}}

\def\weakto{\Rightarrow}

\def\R{{\mathbb R}}

\def\N{{\mathbb N}} 
\def\B{{\mathbb B}}

\def\calD{\mathcal D}
\newcommand{\calA}{{\mathcal A}}
\newcommand{\calB}{{\mathcal B}}

\newcommand{\calF}{{\mathcal F}}
\newcommand{\calG}{{\mathcal G}}

\newcommand{\calL}{{\mathcal L}}
\newcommand{\calM}{\mathcal{M}}
\newcommand{\calN}{\mathcal{N}}
\newcommand{\calP}{{\mathcal P}}

\newcommand{\calW}{\mathcal{W}}

\renewcommand{\d}{{\rm d}}

\newcommand{\ddelta}[1]{\delta_{\pp{#1}}}

\date{}

\def\R{\mathbb{R}}
\def\C{\mathbb{C}}
\def\N{\mathbb{N}}

\def\1{\ind}
\renewcommand{\i}{{\mathrm i}}

\newcommand{\rmIm}{{\mathrm Im}}

\newcommand{\rmsc}{\mathrm{sc}}

\newcommand{\supp}{\mathrm{supp}}

\title{Classical free Brownian motion from the finite free heat flow}
\author{Yizao Wang}
\address{ Department of Mathematical Sciences, University of Cincinnati, 2815 Commons Way, Cincinnati, OH, 45221-0025, USA}
\email{yizao.wang@uc.edu}
\begin{document}\sloppy
\begin{abstract}
A classical free Brownian motion is a classical Markov process that shares the same transition kernels with a free Brownian motion in free probability theory. 
We construct a Markov process that jumps between the evolving roots of a real-rooted polynomial of degree $n$ under the finite free heat flow. Under weak convergence assumptions on the initial data, we prove that this process converges to a classical free Brownian motion as $n\to\infty$. The proof is based on convergence of generators and a martingale problem characterization of classical free Brownian motion.
\end{abstract}
\maketitle
\section{Introduction and main result}
\citet{biane98processes} found an intrinsic connection between non-commutative Markov processes with free increments and a family of classical Markov processes, so that the associated non-commutative and classical processes have the same transition operators. We refer to the classical Markov processes in this correspondence as the {\em classical free Markov processes}.  The most prominent example is the classical free Brownian motion \citep{biane98processes,bozejko97qGaussian} associated to the free Brownian motion \citep{biane97freeBm}. If we let $(\mathbb S_t)_{t\ge0}$ denote a free Brownian motion, then each $\mathbb S_t$ has the semi-circular law with variance $t$,
and $\mathbb S_r$ and $\mathbb S_t - \mathbb S_r$ are freely independent for all $0\le r<t<\infty$. 
Here, we take a more general view following \citet{anshelevich13generators}, who studied the non-commutative stochastic processes with free increments in the form of 
\[
\mathbb X_t:=\mathbb X_0+\mathbb S_t, t\ge0,
\]
 where $(\mathbb S_t)_{t\ge0}$ is as before, $\mathbb X_0$ is a non-commutative random variable with law $\nu$, and $\mathbb X_0$ is freely independent from $(\mathbb S_t)_{t\ge0}$. Biane computed the transition kernels $(\calP^\nu_{s,t})_{0\le s<t<\infty}$ of these processes, which we shall recall below, and established that these kernels satisfy the Chapman--Kolmogorov equation. Anshelevich viewed the process $(\mathbb X_t)_{t\ge0}$ as a free Brownian motion starting from $\mathbb X_0$ with law $\nu$, and in particular he computed explicitly the generators $\calA_t^\nu$ of this process, which will be crucial in our analysis here. Note that for the non-commutative process $(\mathbb X_t)_{t\ge0}$, the initial law $\nu$ has an influence on all the transition operators and generators. This is in contrast to the classical Brownian motion.

In particular, for every probability measure $\nu$, the transition kernels $(\calP^\nu_{s,t})_{0\le s<t<\infty}$ and an arbitrary initial law $\lambda$ determine a classical Markov process. We let $X^{\nu,\lambda}$ denote such a process, and refer to it as a {\em classical free Brownian motion}.  We think of $\nu$ as the environment of the process, in the sense that $\calP_{s,t}^\nu$ depends on $\nu$. A special choice of the initial law is $\lambda = \nu$, and in this case it has been well known that the law of $X^{\nu,\nu}(t)$ is $\nu\boxplus \mu_{\rmsc}^{\boxplus t}$, where $\boxplus$ denotes the free convolution \citep{voiculescu85symmetries,nica06lectures}.

When $\nu =\lambda= \delta_0$, we refer to $X^* = X^{\delta_0,\delta_0}$ as a {\em standard} classical free Brownian motion. This is the most prominent example first studied in the literature \citep{bozejko97qGaussian}. It is a time-inhomogeneous Markov process starting from 0, with semi-circular law $\mu_{\rmsc}^{\boxplus t}$ at time $t>0$, and the transition probability density function given by
\[
\mathsf p_{s,t}(x,y) = \frac{t-s}{2\pi}\frac{\sqrt{4t-y^2}}{(t-s)^2+tx^2-(s+t)xy+sy^2}\inddd{|y|<2\sqrt t}, \mbox{ for } x\in(-2\sqrt s,2\sqrt s).
\]
The process $X^*$ and a few other classical free Markov processes, referred to as the classical free Meixner processes in the literature, also appeared independently in the studies of so-called quadratic harnesses by \citet{bryc05conditional}. It turned out that these processes can be put in an even larger family of so-called Askey--Wilson Markov processes introduced by the same two authors \citep{bryc10askey}, and in the past decade they have played a crucial role in integrable probability \citep{bryc19limit,corwin24stationary}.

The introduction of these classical free Markov processes is functional analytic in nature. We are instead interested in probabilistic interpretations of these processes, and we are aware of very few results. In particular, we are interested in stochastic models in which classical free Markov processes arise in the limit, and to the best of our knowledge such a limit theorem has not been completely explicit in the literature. 
The main goal of this paper is to explain how the classical free Brownian motion arises as the limiting process from the finite free heat flow.
\subsection{The finite free heat flow}
We call the dynamics induced by the following semigroup, for each $n\in\N$, the finite free heat flow:
\[
T_t\topp n :=\exp\pp{-\frac t{2n}\pp{\frac{\d^2}{\d x^2}}}, t\ge 0.
\]
Then, given a polynomial $p_{n,0}(x)$ of degree $n$, the finite free heat flow induces a family of polynomials
\[
p_{n,t}(x) = T_t\topp np_{n,0}(x),  t\ge 0,
\]
where for each $t>0$ the polynomial $p_{n,t}$ remains  a polynomial of degree $n$. Most importantly, if $p_{n,0}$ has only real roots, then $p_{n,t}$ also has only real roots. 

The name of the dynamics reflects the connection to the finite free probability theory \citep{marcus22finite}. Let $p_{n,0}$ and $p_{n,t}$ be as above, and $\boxplus_n$ be the finite free (additive) convolution of polynomials of degree $n$.  Then, we also have
\[
p_{n,t} = p_{n,0}\boxplus_n H_{n,t},
\]
where $H_{n,t}(x):=T_t\topp n x^n$ is the rescaled monic Hermite polynomial.

We are interested in the evolution of the real roots of a polynomial under the finite free heat flow. 
That is, writing $\vec a_n = (a_{n,1},\dots,a_{n,n})$ and explicitly
\[
p_{n,0}(x)\equiv p_{\vec a_n}(x) := \prodd i1n (x-a_{n,i}),
\]
we are interested in the roots of $p_{\vec a_n,t} = p_{\vec a_n}\boxplus_n H_{n,t}$, denoted by 
\[
\vec r_n(t) = (r_{n,1}(t),\dots,r_{n,n}(t)), t\ge0.
\] Note that this vector is determined by the dynamics and the initial roots $\vec a_n$, and we assume throughout that the roots are listed in the increasing order $r_{n,1}(t)<\cdots<r_{n,n}(t)$. We assume all the initial roots to be distinct for now. We shall refer to $(\vec r_n(t))_{t\ge0}$ as the {\em finite free root flow}.  They satisfy in particular the following ODE system
\equh\label{eq:frozen DBm}
r_{n,i}'(t) = \frac1n\sum_{j=1,\dots,n, j\ne i}\frac1{r_{n,i}(t)-r_{n,j}(t)}, \quad i=1,\dots,n, t>0.
\eque
The proof can be found in \citet[Lemma 2.4]{csordas94lehmer}. Recent discussions of this fact can be found in \citet{garza-vargas26finite} and \citet{voit25differential}. 
Let $\mu[p]$ denote the empirical root measure of the polynomial $p$. The following results have been known regarding the empirical root measures: if 
\[
\mu [p_{\vec a_n}] = \frac1n\summ i1n \delta_{a_{n,i}}\weakto\nu, 
\]
then
\[
 \mu\bb{p_{\vec a_n,t}}\equiv \mu\bb{p_{\vec a_n}\boxplus_n H_{n,t}} = \frac1n\summ i1n \delta_{r_{n,i}(t)}\weakto \mu_t^\nu:=\nu\boxplus \mu_\rmsc^{\boxplus t},
\]
as $n\to\infty$, where  
\[
\mu_\rmsc(\d x) = \frac{\sqrt{4-x^2}}{2\pi}\inddd{x\in(-2,2)}\d x
\]
is the standard semi-circular law. 
For the convergence above, see \citet[Theorem~2.13 and Remark~2.15]{kabluchko25lee}, \citet[Theorem~1.3]{fujie26regularity}, and \citet[Theorem~1.1]{voit22limit}. The same finite free polynomial flow has been studied with random initial polynomials by \citet{hall25zeros}, and with initial polynomials in the form of $P(x)^{n/\mathrm{deg}(P)}$ (with $P$ fixed) by \citet{hofert26zeros}. The same limiting flow $(\mu_t^\nu)_{t\ge0}$ also arises in the large-dimensional limits of empirical spectral measure processes, after deterministic changes of time and space \citep{chan92wigner,rogers93interacting}. For related work on rectangular finite free heat flow, see \citet{cuenca26rectangular}. 
\subsection{Main result}\label{sec:main result}

We introduce a Markov process that at each time $t$ takes values in $\{r_{n,1}(t),\dots,r_{n,n}(t)\}$. We first define a Markov process taking values in $\{1,\dots,n\}$. Set, for $\vec r = (r_1,\dots,r_n)$ with all $n$ coordinates having distinct values,
\[
Q\topp n_{i,j}(\vec r ) := \begin{cases}
\displaystyle\frac1n\frac1{(r_{i}-r_{j})^2}, & \mbox{ if } i\ne j,\\\\
\displaystyle-\frac1n\sum _{k\ne i}\frac1{(r_{i}-r_{k})^2}, & \mbox{ if } i=j.
\end{cases}
\]
We let $I_n = (I_n(t))_{t\ge 0}$ denote the Markov process on $\{1,\dots,n\}$ determined by the generators $(\calL\topp n_t)_{t\ge0}$ given by 
\equh\label{eq:I generator}
\bb{\calL\topp n_{t}f}(i) =\sum_{j\ne i}Q\topp n_{i,j}(\vec r_n(t))(f(j)-f(i)) =  \frac1n\sum_{j\ne i}\frac{f(j)-f(i)}{(r_{n,j}(t)-r_{n,i}(t))^2}, i=1,\dots,n,
\eque
and an initial law to be prescribed later. We refer to the process $I_n$ as the {\em rank} process. We are mostly interested in the {\em tagged-root} process 
\[
X_n(t):=r_{n,I_n(t)}(t), t\ge 0.
\]
In words, the process $X_n\equiv (X_n(t))_{t\ge 0}$ records the trajectory of a particle jumping between $n$ different trajectories of the roots of the polynomials, where the jumps are encoded by the rank process $I_n$.

The main result of the paper is the next theorem on the convergence to a classical free Brownian motion $X^{\nu,\lambda}$ introduced earlier. Let $\Phi_t^{\nu}(x) := \proba(X^{\nu,\nu}(t)\le x)$. We emphasize that $\Phi_t^\nu$ does not depend on $\lambda$.
\begin{theorem}\label{thm:2}
Consider a finite free root flow $(\vec r_n(t))_{t\ge0}$ with the initial roots satisfying
$r_{n,1}(0)<\cdots<r_{n,n}(0)$ for all $n\in\N$ and 
\equh\label{eq:cond nu}
\nu_n:=\frac1n\sum_{i=1}^n\delta_{r_{n,i}(0)}\weakto\nu,
\eque
as $n\to\infty$, where $\nu$ is a probability measure. Consider a family of rank processes $(I_n)_{n\in\N}$ described above with initial distributions satisfying 
\equh\label{eq:cond lambda}
X_n(0) = r_{n,I_n(0)}(0)\weakto\lambda,
\eque
as $n\to\infty$.
Then, for every $\epsilon>0$, 
\equh\label{eq:joint}
\pp{\pp{\frac1nI_n(t)}_{t\ge \epsilon},\pp{X_n(t)}_{t\ge0}}\weakto\pp{\pp{\Phi^\nu_t\pp{X^{\nu,\lambda}(t)}}_{t\ge\epsilon},\pp{X^{\nu,\lambda}(t)}_{t\ge0}}
\eque
in $D([\epsilon,\infty))\times D([0,\infty))$, equipped with the Skorokhod $J_1$ topology. If in addition $\nu$ has no atoms, then $\epsilon$ in \eqref{eq:joint} can be taken to be zero. 
\end{theorem}

The restriction of the rank process to $[\epsilon,\infty)$ for $\epsilon>0$ cannot be removed without additional assumptions. See Remark \ref{rem:zero}. 
\begin{remark}
The most natural choice of the initial law of the rank process $I_n(0)$ is the uniform distribution over $\{1,\dots,n\}$, and in this case necessarily $\lambda = \nu$. 
However, in general for any probability measures $\nu,\lambda$, one may take certain initial roots $(\vec r_n(0))_{n\in\N}$ and the initial laws of $(I_n(0))_{n\in\N}$ such that \eqref{eq:cond nu} and \eqref{eq:cond lambda} both hold. In particular, $\lambda$ does not have to be absolutely continuous with respect to $\nu$.  A remarkable situation is when $\lambda = \delta_a$ for some $a\notin\supp(\nu)$, and in this case the paths of $X^{\nu,\delta_a}$ start from outside $\supp(\nu)$ and may remain outside $\supp (\nu\boxplus \mu_\rmsc^{\boxplus t})$ for a positive time, while the paths of $X^{\nu,\nu}$ remain within $\supp(\nu\boxplus\mu_\rmsc^{\boxplus t})$ from the beginning. See end of Section \ref{sec:classical free Bm} for more discussions. 
\end{remark}
In Theorem \ref{thm:2}, we assume the roots are all distinct at time zero. This assumption can be further relaxed to allow some of the initial roots to be identical. We do not pursue a fully general statement, but single out a special case. That is, the initial polynomial is 
\[
p_{n,0}(x) = x^n.
\] 
We refer to the root flow in this example as the {\em Hermite root flow}: in this case, at any time $t>0$, $p_{n,t} = (t/n)^{n/2} H_n(\sqrt{n/t}x)$ where $H_n$ is the probabilist's Hermite polynomial of degree $n$.  This time, $I_n(0)$ cannot be properly defined as the process has infinitely many jumps in the neighborhood of the time origin. However, the tagged-root process for the Hermite root flow, denoted by $X_n^*$, is defined by an approximation argument near the time origin, and we can still show that $X_n^*\weakto X^{\delta_0,\delta_0}$. Section \ref{sec:hermite entrance} is devoted to this example.

The proof of the convergence $X_n\weakto X^{\nu,\lambda}$ in $D([0,\infty))$ in Theorem \ref{thm:2} is based on a martingale problem, which we first develop in Section \ref{sec:martingale}. Then, to prove the convergence we essentially need to control the difference between, for a family of test functions $f$,  the generator of the process $X^{\nu,\lambda}$ due to \citet{anshelevich13generators}, which is given by
\[
\calA_t^\nu f(x) = \partial_x\int\frac{f(x)-f(y)}{x-y}\pp{\nu\boxplus\mu_\rmsc^{\boxplus t}}(\d y) = \int K_f(x,y)\pp{\nu\boxplus \mu_\rmsc^{\boxplus t}}(\d y), 
\]
with a kernel function $K_f$ defined later in \eqref{eq:K_f}, and the generator of the process $X_n$, given by
\[
\calA_{n,t}f(r_{n,i}(t)) = \int K_f(r_{n,i}(t),y)\mu_{n,t}(\d y)  - \frac1{2n}f''(r_{n,i}(t)),
\]
where $\mu_{n,t} = n\inv\summ i1n \delta_{r_{n,i}(t)}\weakto \nu\boxplus\mu_\rmsc^{\boxplus t}$ (see Lemma \ref{lem:generator}). 
We also establish the tightness of the sequence $(X_n)_{n\in\N}$. The proof of the joint convergence of $(I_n/n,X_n)$ then follows by a continuous mapping argument. The proofs of limit theorems are provided in Section \ref{sec:proofs}.

\subsection{Related work from random matrix theory}\label{sec:RMT}
The finite-dimensional labeled dynamics underlying $X_n$ can be related to frozen Dunkl processes and has already been studied. Let $(\vec Z_n^*(t))_{t\geq 0}$ be the frozen type-$A_{n-1}$ Dunkl process of \citet{voit22limit}, and set $\vec Z_n=n^{-1/2}\vec Z^*_n$ to match our spatial normalization. 
If $I_n(0)$ is uniformly distributed over $\{1,\dots,n\}$ and $\vec Z_n(0)$ is a uniformly random permutation of $\vec r_n(0)$, then almost surely,
\[
  \{Z_{n,1}(t),\ldots,Z_{n,n}(t)\}  =  \{r_{n,1}(t),\ldots,r_{n,n}(t)\},
  \qquad t\geq 0,
\]
and, for every $i=1,\ldots,n$,
\[
  \bigl(Z_{n,i}(t)\bigr)_{t\geq 0}  \stackrel{\mathrm d}{=}  \bigl(X_n(t)\bigr)_{t\geq 0}.
\]
Thus, $X_n$ may be viewed as a tagged coordinate of a spatially rescaled frozen type-$A$ Dunkl process.  Voit and Woerner study large-$n$ limits at the level of empirical root measures, at which the random permutation dynamics is invisible.  They do not consider the tagged-particle process.  
In the case of the Hermite root flow, the permutation jump dynamics in the freezing regime was studied by \citet{andraus20dunkl}.
Projecting its generator onto the rank of a fixed particle yields our rank generator; see \citet[Section~6, equations~(101)--(102)]{andraus20dunkl}.

There are a few other closely related models in the literature. A dynamics intrinsically connected to our model is 
\equh\label{eq:B}
\frac{\d \vec z}{\d t} = Q(\vec r(t))\vec z.
\eque
(We omit the dependence on $n\in\N$ for the rest of this discussion.)
In fact, the transposed transition matrices $(P_{s,t}\topp n)^T$ of the rank process $I_n$ satisfy the above ODE (see \eqref{eq:PQ} below). The matrix $(P_{s,t}\topp n)^T$ can be interpreted as the Jacobian of the changes of the finite free root flow at time $t$ caused by perturbation at time $s$, and this matrix has shown up in the recent investigation of \citet{garza-vargas26finite}; see Section \ref{sec:finite free} and Remark \ref{rem:GVSS} for more details. \bigskip

Moreover, the finite free heat flow $(\vec r_n(t))_{t\ge0}$ is intrinsically related to the $n$-dimensional Dyson Brownian motion. Recall the ODE satisfied by the finite free root flow in \eqref{eq:frozen DBm}. For the convenience of comparison, we simply write 
\equh\label{eq:A}
\frac{\d\vec r}{\d t} = \vec b(\vec r).
\eque
It is well-known that the SDE corresponding to \eqref{eq:A}, 
\[
\d\vec \lambda = \vec b(\vec\lambda)\d t + \sigma \d\vec \B_t,
\]
 characterizes the $n$-dimensional Dyson Brownian motion. 
Consider
\equh\label{eq:E}
\frac{\d \vec u}{\d t} = Q\pp{\vec \lambda(t)} \vec u.
\eque
This is the same ODE as in \eqref{eq:B} except that now $Q$ is evaluated at the random eigenvalues $\vec\lambda(t)$ of the Dyson Brownian motion instead of the deterministic finite free root flow $\vec r(t)$. 
This dynamics showed up in the investigation of the eigenvector moment flow of \citet{bourgade17eigenvector}. Consider a Hermitian matrix Brownian motion $H_t=A+W_t$, with $A$ having simple spectrum, and write $\lambda_i(t)$ and $v_i(t)$ for its ordered eigenvalues and corresponding normalized eigenvectors. They studied a certain dynamics of $n$ particles on the $N$ eigenvalue ranks \citep[Theorem~3.1]{bourgade17eigenvector}. The dynamics of a single particle corresponds to our model here.

For the same model, fixing an initial eigenvector $v_j$, \citet{allez14eigenvectors} considered the overlaps of eigenvectors
\[
O_{i\mid j}(t)
:=
\mathbb E\pp{
\left|\left\langle v_i(t),v_j(0)\right\rangle\right|^2\mmid\vec\lambda(s), s\in[0,t]}.
\]
Their overlap equation is
\[
\partial_tO_{i\mid j}(t)
=
\sum_{k\ne i}Q_{ik}(\vec\lambda(t))\pp{O_{k\mid j}(t)-O_{i\mid j}(t)},
\qquad O_{i\mid j}(0)=\delta_{ij}.
\]
This is precisely \eqref{eq:E}. In fact, as they pointed out in the paper, this overlap evolution equation in fact already appeared in \citet[equation~(4.7)]{wilkinson95brownian}.
In particular, \citet{allez14eigenvectors} analyzed the associated local spectral measures
\[
\sum_{i=1}^n O_{i\mid j}(t)\,\delta_{\lambda_i(t)}.
\]
If we define the rank process $I_n$ accordingly,  this measure is exactly the conditional one-time distribution of the tagged eigenvalue $\lambda_{I_n(t)}(t)$ started from rank $j$. 
Under their assumptions on the initial matrices, whose empirical spectral measures satisfy
$n\inv\sum_{i=1}^n\delta_{\lambda_i(0)}\weakto\nu$,
they showed that if the initial eigenvalue $\lambda_j(0)\equiv \lambda_{j_n}(0)\to a$ for a sequence $(j_n)_{n\in\N}$, then  the limit of the measure above is the corresponding transition law of classical free Brownian motion $\calP_{0,t}^\nu(a,\d y)$ in our notation
\citep[Section~3.2, equation~(3.15) and the following paragraph]{allez14eigenvectors}. 
This is not yet a result on the process convergence to the classical free Brownian motion, but it strongly suggests so. As noted by \citet{allez14eigenvectors}, the limiting local resolvent in their equation~(3.15) had already been obtained by \citet{shlyakhtenko96random} using free probability. 

The ODE \eqref{eq:E} also shows up in the work of  \citet{landon19fixed}, when studying an initial-data interpolation for two Dyson Brownian motion trajectories with initial data $\vec a$ and $\vec b$, respectively. 
\bigskip

Another closely related model showed up in the studies of \citet{bourgade16fixed}. Let $(x_{n,i}(t))_{i=1,\dots,n}$ and $(y_{n,i}(t))_{i=1,\dots,n}$ be two Dyson Brownian motions (in fact, the corresponding Ornstein--Uhlenbeck processes) driven by the same Brownian motions, and consider a rate matrix $\calB\topp n(t)$  depending on two solutions of Dyson Brownian motion. More precisely, the $(i,j)$-entry of the matrix is approximately 
\[
\calB_{i,j}
\topp n(t) \approx\frac1n\frac1{(x_{n,i}(t)-x_{n,j}(t))(y_{n,i}(t)-y_{n,j}(t))}, i\ne j,
\]
and $\calB_{i,i}\topp n(t) = -\sum_{j\ne i}\calB_{i,j}\topp n(t)$.
They considered the matrix $U_n(s,t)$ satisfying
\[
\partial_tU_n(s,t) = \calB\topp n(t)U_n(s,t) \qmand U_n(s,s) = \mathbb I_n.
\] 
(Compare with our \eqref{eq:variational} below.)
Letting $\gamma_{n,j}$ be such that $\mu_{\mathrm {sc}}([-2,\gamma_{n,j}]) = (j-1/2)/n, j=1,\dots,n$, their result essentially implies that with $(i_n)_{n\in\N}$ such that $\limn\gamma_{n,i_n} = x\in(-2,2)$, 
\[
\summ j1n \bb{U_n(s,t)}_{j,i_n}\delta_{\gamma_{n,j}} \to \mathsf p^{\mathrm{free OU}}_{t-s}(x,y)\mu_{\mathrm{sc}}(\d y),
\]
in probability as $n\to\infty$, where $\mathsf p^{\mathrm{free OU}}_t(x,y)$ is the transition probability density function of the classical free Ornstein--Uhlenbeck process with respect to $\mu_\rmsc$
 \citep[Theorems 3.3 and 3.4]{bourgade16fixed}
 (they established a much more precise estimate of this convergence). They did recognize the connection to classical free Markov processes \citep[Remark 1]{bourgade16fixed}, although this connection is not the main focus of their paper. The above convergence is essentially of the same type as our convergence of the tagged root process at the level of transition kernels.
\bigskip

The same matrix $Q(\vec r(t))$ also appears in the linear SDE
\[
\d \vec\zeta = Q(\vec r(t))\vec\zeta\d t + \sigma \d \vec\B_t.
\]
(Compare with \eqref{eq:B}.)
This SDE describes the Gaussian fluctuation process of Dyson Brownian motion in the freezing regime, which arises as the limit of appropriately normalized Dyson Brownian motion as the inverse temperature $\beta\uparrow\infty$, studied first by \citet{voit21functional} (who assumed distinct initial values) and  \citet{gorin24universal} for the case of the Hermite root flow. See also \citet{andraus19central} for an earlier development on the marginal convergence.

\subsection*{The paper is organized as follows.} Section \ref{sec:preliminary} provides background on the finite free heat flow and the classical free Markov process. The martingale problem for the classical free Brownian motion is proved in Section \ref{sec:martingale}. The limit theorems are proved in Section \ref{sec:proofs}. 

\subsection*{Statement on the use of AI} The author had been looking for a stochastic model that scales to the classical free Brownian motion since his work on the topic \citep{bryc16local}, but without success. When he switched to ChatGPT for assistance in summer 2026, ChatGPT first found the key references of \citet{allez14eigenvectors}, \citet{bourgade16fixed}, and \citet{bourgade17eigenvector}, whose results strongly suggested that such a limit theorem can be established from dynamics based on Dyson Brownian motion (as explained in Section \ref{sec:RMT}, the transition kernels have been explicitly computed in a couple of cases). The author decided not to pursue this direction. The author then asked ChatGPT to search the literature on finite free probability, and ChatGPT found another key reference of \citet{garza-vargas26finite}. The reference itself does not focus on the Markov-process point of view, and the connection to the original question was pointed out by ChatGPT (see Remark \ref{rem:GVSS}). ChatGPT played the role of co-author in the preparation of this paper. The author wrote the paper and takes full responsibility for the content. 

\subsection*{Acknowledgement} YW would like to thank W{\l}odek Bryc for several discussions, and in particular for suggesting a search for connections to finite free probability in summer 2026.  YW's research was partially supported by the Simons Foundation (MP-TSM-00002359).

\section{Preliminary background}\label{sec:preliminary}
We first explain how the generators of the rank process $I_n$ can be related to the Jacobians of the ODE \eqref{eq:frozen DBm} in Section \ref{sec:finite free}. We then provide some background on the classical free Markov process in Section \ref{sec:classical free Bm}.
\subsection{The finite free root flow}\label{sec:finite free}
In Section \ref{sec:main result}, we define the Markov process through its generator $Q(\vec r_n(t))$. The construction is guaranteed by the fact that $t\mapsto Q(\vec r_n(t))$ is continuous and bounded on compact time intervals. 
Here, we present an alternative approach. We identify the transition matrices of the rank process $I_n$ with the transposes of the Jacobians of the finite free root flow.

For $t\ge 0$, let
\[
\vec\lambda_{n,t}(\vec a_n) = \pp{\lambda_{n,1,t}(\vec a_n),\dots,\lambda_{n,n,t}(\vec a_n)}
\]
denote the vector of increasing roots of $T_t\topp n p_{\vec a_n}$. We have
$\vec\lambda_{n,s}\spp{\vec\lambda_{n,t}(\vec a_n)} = \vec\lambda_{n,s+t}(\vec a_n)$.
Given a finite free root flow $(\vec r_n(t))_{t\ge0}$, one may consider the associated Jacobian with respect to the changes at time $t$ caused by the perturbation at time $s$:
\[
J_{s,t}\topp n:= \calD\vec\lambda_{n,t-s}(\vec r_n(s)), 0\le s<t.
\]
That is, 
\[
J_{s,t}\topp n(i,j) = \frac{\partial \lambda_{n,i,t-s}(\vec x)}{\partial x_j}\bigg\vert_{\vec x = \vec r_n(s)}\equiv  \frac{\partial r_{n,i}(t)}{\partial r_{n,j}(s)}, \quad 1\le i,j\le n, s\le t. 
\]
The Jacobian $J_{s,t}\topp n$ in fact solves the variational equation of the ODE of the finite free root flow \eqref{eq:frozen DBm}. 
Indeed, rewriting \eqref{eq:frozen DBm} as 
\[
r_{n,i}'(t) = b_{n,i}(\vec r_n(t)) \qmwith b_{n,i}(\vec r_n(t)) = \frac1n\sum_{j=1,\dots,n, j\ne i}\frac1{r_{n,i}(t)-r_{n,j}(t)}, \quad i=1,\dots,n, t>0,
\]
and $\vec b_n(\vec r_n(t)) = (b_{n,i}(\vec r_n(t)))_{i=1,\dots,n}$, we have 
 \[
 Q\topp n(\vec r_n(t)) = \calD \vec b_n(\vec r_n(t)),
 \]
and 
\equh\label{eq:variational}
\partial_t J_{s,t}\topp n =  Q\topp n(\vec r_n(t))J\topp n_{s,t} \qmand J\topp n_{s,s} = \mathbb I_n, 0\le s<t.
\eque
Now, introduce
\[
P\topp n_{s,t}:=\pp{J_{s,t}\topp n}^T, 0\le s<t.
\]
The next proposition shows that these matrices are precisely the transition matrices of the rank process $I_n$ introduced in Section \ref{sec:main result}.

\begin{proposition}\label{prop:1}
Assume that $r_{n,1}(0)<\cdots<r_{n,n}(0)$. Then, the matrices $(P\topp n_{s,t})_{0\le s\le t}$ are transition matrices of the rank process $I_n$, whose instantaneous generator at time $t$ is $Q\topp n(\vec r_n(t))$. That is,
\equh\label{eq:PQ}
\partial_tP_{s,t}\topp n = P_{s,t}\topp nQ\topp n\pp{\vec r_n(t)}, \quad P_{s,s}\topp n = \mathbb I_n. 
\eque
\end{proposition}
The proof of Proposition \ref{prop:1} is straightforward. 
 Since $Q\topp n$ is symmetric and $Q\topp n\ind  =\ind ^TQ\topp n = 0$, it follows that $J\topp n_{s,t}\ind =\ind$ and $\ind^T J\topp n_{s,t} = \ind^T$. Since the off-diagonal entries of $Q\topp n$ are all non-negative, it follows that all the entries of $J\topp n_{s,t}$ are non-negative. That is, $J_{s,t}\topp n$ is doubly stochastic. The fact that $J_{s,t}\topp n$ is doubly stochastic has been actually proved in a more general setup by \citet{garza-vargas26finite}. It remains to check the Chapman--Kolmogorov equation, which is straightforward as $J_{r,t}\topp n = J_{s,t}\topp nJ_{r,s}\topp n$. 

\begin{remark}\label{rem:GVSS} The fact that $J_{s,t}\topp n$ for a fixed pair of $s,t$ is doubly stochastic was  already observed by  \citet[Lemma 3.4]{garza-vargas26finite} in a much more general setting.  In their notation, our $J_{s,t}\topp n$ is their $E_\beta[\alpha\mid\gamma]$,
where $\alpha=\vec r_n(s)$, $\beta$ is the vector of increasing roots of $T_{t-s}\topp n x^n$, and $\gamma=\Omega_{\boxplus_n}(\alpha,\beta)=\vec r_n(t)$. 
  They considered the finite free convolution of more general polynomials. For the finite free heat flow, the ODE \eqref{eq:frozen DBm} plays a key role in their analysis  \citep[Lemma 2.5]{garza-vargas26finite}.  
They did not take the Markov process point of view of $J_{s,t}\topp n$. Instead, they pointed out the doubly stochastic structure of the Jacobian, and mentioned that this could be interpreted as a coupling of probability distributions \citep[Section 3.3]{garza-vargas26finite} (their discussions are not restricted to the finite free heat flow). In our notation,  they are concerned with
\[
\Lambda_{s,t}\topp n:=\frac1n\summ i1n \summ j1n J_{s,t}\topp n(i,j)\ddelta{r_{n,j}(s), r_{n,i}(t)},
\]which is in fact the bivariate distribution of the tagged root process at times $s$ and $t$, if started with the uniform distribution over $\{r_{n,1}(s),\dots,r_{n,n}(s)\}$ at time $s$. 
\end{remark}

\subsection{Classical free Markov processes}\label{sec:classical free Bm}
We shall work with the classical free Brownian motion, a prominent example of a classical free Markov process  in the framework of \citet{biane98processes}.
This is a family of time-inhomogeneous Markov processes denoted by $X^{\nu,\lambda} = (X^{\nu,\lambda}_t)_{t\ge 0}$, where $\nu,\lambda$ are probability measures on $\R$. We let $\calP^\nu = (\calP_{s,t}^\nu)_{0\le s<t<\infty}$ denote the transition kernels and $\lambda$ the initial law of the process.

The transition kernels $\calP^\nu$ originate in the work of
\citet{biane97free,biane98processes}. The process naturally associated
with Biane's construction has initial law $\nu$, namely
$X^{\nu,\nu}$. Here we use the same transition kernels with a more
general initial law $\lambda$.
Later on, \citet{anshelevich13generators} studied its generator $\calA^\nu_t$, and showed that
\[
\calA_t^\nu f(x) = \partial_x\int\frac{f(x)-f(y)}{x-y}\pp{\nu\boxplus\mu_\rmsc^{\boxplus t}}(\d y), t> 0,
\]
for all $f\in C_c^\infty(\R)$. 
 For the proofs of our main theorem, we do not need an explicit expression of $\calP^\nu$ but shall only use the generator. 
(In fact, it is worth noting that the derivation of the expression of $\calA_t^\nu$ does not involve the explicit expression of $\calP_{s,t}^\nu$ either.) Here, $\boxplus$ is the free convolution and we shall write
\[
\mu_t\equiv \mu_t^\nu :=\nu\boxplus \mu_\rmsc^{\boxplus t}, t\ge 0,
\]
in short. 
We regard $\calA_t^\nu$ as the generator on $C_0(\mathbb R)$. Below, we use its restriction to the test space $\calW$ introduced in \eqref{eq:calW}. 
Thus,
\[
\calW \subseteq D(\calA_t^\nu),
\qquad
\calA_t^\nu(\calW)\subseteq C_0(\mathbb R).
\]
See \citet[Theorem 15]{anshelevich13generators}.

We provide a brief description of $(\calP^\nu_{s,t})_{0\le s<t}$ here. 
First, it is well-known that for every $\nu$, there exists a noncommutative probability space $(\calM,\tau)$, a free Brownian motion $(\mathbb S_t)_{t\ge 0}$ and a noncommutative random variable $\mathbb X_0$ with law $\nu$ and freely independent from the free Brownian motion. Then, set
\[
\mathbb X_t:=\mathbb X_0+\mathbb S_t, t>0.
\]
We have that 
\[
\calL(\mathbb X_t) = \mu_t^\nu = \nu\boxplus\mu_\rmsc^{\boxplus t}, t\ge 0.
\]
Biane first studied the law $\mu_t^\nu$ \citep{biane97free}, and then  proved that the so-constructed random variables $(\mathbb X_t)_{t\ge 0}$ have the noncommutative Markov property
 \citep{biane98processes}.  In particular, he discovered the key identity that
\equh\label{eq:NC markov}
\esp_\tau\pp{\frac1{z-\mathbb X_t}\mmid \calN_s} = \frac1{F_{s,t}^\nu(z)-\mathbb X_s}, \quad 0\le s<t, z\in\C^+,
\eque
for an analytic function $F_{s,t}^\nu:\C^+\to\C^+$. Here, $\calN_s$ is the von Neumann algebra generated by $(\mathbb X_r)_{r\in[0,s]}$, and $\esp_\tau(\cdot\mid\calN_s)$ is the $\tau$-preserving conditional expectation onto $\calN_s$.
The function $F_{s,t}^\nu$ is the unique analytic function from $\C^+$ to $\C^+$ satisfying, with $G_\mu(z) = \int (z-x)\inv \mu(\d x)$ being the Cauchy transform of $\mu$, 
\[
G_{\mu_t^\nu}(z)
   =G_{\mu_s^\nu}\pp{F_{s,t}^\nu(z)},
\qmand
\lim_{y\to\infty}\frac{F_{s,t}^\nu(\i y)}{\i y}=1,
\]
and  turns out to be
\equh\label{eq:F}
F_{s,t}^\nu(z) = z-(t-s)G_{\mu_t^\nu}(z).
\eque
Then, Biane recognized that
\[
\frac1{F_{s,t}^\nu(z)-x}
\]
for every $x\in\R$
is a Cauchy transform (as a function of $z\in\C^+$) of a probability measure by the Herglotz--Nevanlinna Theorem. We denote this probability measure by  $\calP_{s,t}^\nu(x,\cdot)$. That is,
\equh\label{eq:resolvent F}
\int \frac1{z-y} \calP_{s,t}^\nu(x,\d y) = \frac1{F_{s,t}^\nu(z)-x}, \mfa x\in\R, z\in\C^+.
\eque
Write $\calP_{s,t}^\nu f(x) = \int f(y)\calP_{s,t}^\nu(x,\d y)$. Biane also proved that
\[
\calP_{s,t}^\nu(C_0(\R))\subset C_0(\R),
\]
and that  the Chapman--Kolmogorov equation holds:
\[
\calP_{r,s}^\nu\calP_{s,t}^\nu = \calP_{r,t}^\nu, \mfa 0\le r<s<t.
\]
So this is a family of Markov transition kernels. 

From now on, we let $X^{\nu,\lambda}$ denote a Markov process with transition kernels $(\calP^\nu_{s,t})_{0\le s<t}$ and initial law $\lambda$. Set $\calF_s^{\nu,\lambda} = \sigma(X^{\nu,\lambda}_r: r\in[0,s])$. We have
\[
\esp\pp{\frac1{z-X_t^{\nu,\lambda}}\mmid\calF_s^{\nu,\lambda}} = \frac 1{F_{s,t}^\nu(z)-X_s^{\nu,\lambda}}, \quad \mbox{ almost surely, } 0\le s<t, z\in\C^+,
\]
which is a classical conditional expectation relation parallel to \eqref{eq:NC markov}.  The so-defined Markov process has a version in $D([0,\infty))$. 
\begin{lemma}\label{lem:D}
For every pair $(\nu,\lambda)$, $X^{\nu,\lambda}$  admits a c\`adl\`ag realization (i.e., a Markov process with sample paths in $D([0,\infty))$) with transition kernels $\calP^\nu$ and initial law $\lambda$. 
\end{lemma}
\begin{proof}
By \citet[Proposition~7]{anshelevich13generators}, $\calP_{s,t}^\nu$ is strongly continuous. 
Consider the associated space--time semigroup on
$E=[0,\infty)\times\mathbb R$,
\[
(T_h\varphi)(s,x)
:=
\calP_{s,s+h}^\nu\bigl(\varphi(s+h,\cdot)\bigr)(x).
\]
It then follows that $T_h$ is a strongly continuous conservative Feller semigroup on $C_0(E)$; see also \citet[Remark~12]{anshelevich13generators}. By the standard realization theorem for Feller semigroups \citep[Chapter~4, Theorem~2.7]{ethier86markov},
it admits a c\`adl\`ag realization for every initial distribution. Now, consider the space--time process with initial law $\delta_0\otimes\lambda$ that is a c\`adl\`ag realization, and notice that the second component of this process has the law of the process $X^{\nu,\lambda}$, and again is c\`adl\`ag. This second component then provides  a desired realization. 
\end{proof}

We emphasize that for the classical free Brownian motion, the initial law $\lambda$ can be an arbitrary probability measure on $\R$. The choice of $\lambda = \nu$ is special as then the marginal law of the process $X^{\nu,\nu}$ becomes $\mu_t^\nu$. Indeed, we have $\mu_s^\nu\calP_{s,t}^\nu = \mu_t^\nu$, since
\begin{align*}
G_{\mu_s^\nu\calP_{s,t}^\nu}(z) &= \int_\R\int_\R\frac1{z-y}\calP_{s,t}^\nu(x,\d y)\mu_s^\nu(\d x) = \int_\R\frac1{F_{s,t}^\nu(z)-x}\mu_s^\nu(\d x) \\
&= G_{\mu_s^\nu}(F_{s,t}^\nu(z)) = G_{\mu_t^\nu}(z).
\end{align*}

More generally,
\[
\calL(X_t^{\nu,\lambda})
   =\lambda\calP_{0,t}^\nu,
\]
which in general is not the same as $\mu_t^\nu$. At the same time, if $\supp\lambda\subseteq\supp\nu$, it follows that $\supp\calL(X_t^{\nu,\lambda})   \subseteq\supp(\mu_t^\nu)$. Thus, if $X^{\nu,\lambda}$ starts from within the support of $\nu$, then it remains within the support of $\mu_t^\nu$ at every $t>0$. However, the process may start from outside the support of $\nu$ at time $0$, and remain outside until a later time $t$. This phenomenon has been known. We collect a few explicit formulas in the next remark.

\begin{remark}
A point outside the support of the free measure flow can remain outside for a positive time. Fix $s\ge0$ and $x\notin\supp(\mu_s^\nu)$, and set
\[
\tau_s(x):=\left(\int_\R\frac{\mu_s^\nu(\d y)}{(x-y)^2}\right)^{-1}
\qmand
\psi_{s,t}(x):=x+(t-s)G_{\mu_s^\nu}(x),\qquad t\ge s.
\]
The Cauchy transforms of the free measure flow satisfy the complex Burgers equation
\[
\partial_t G_{\mu_t^\nu}(z)
+G_{\mu_t^\nu}(z)\partial_z G_{\mu_t^\nu}(z)=0,
\qquad t>0,\ z\in\C^+;
\]
see \citet{biane97free}. 
The real characteristic curve starting from $(s,x)$ is
\[
z(t)=x+(t-s)G_{\mu_s^\nu}(x)=\psi_{s,t}(x),
\qquad t\in[s,s+\tau_s(x)],
\]
and this curve is a solution of $z'(t)=G_{\mu_t^\nu}(z(t))$. Along such a curve, the Burgers equation gives $\frac{\d}{\d t}G_{\mu_t^\nu}(z(t))=0$, so its velocity is constant.  The quantity $s+\tau_s(x)$ is the first time at which the characteristic curve meets the support of the free measure flow. Indeed, applying \citet[Lemma~4 and Corollary~3]{biane97free} with initial measure $\mu_s^\nu$, we obtain
$\psi_{s,t}(x)\notin\supp(\mu_t^\nu)$ for all $t\in(s,s+\tau_s(x))$ and $\psi_{s,s+\tau_s(x)}(x)\in\supp(\mu_{s+\tau_s(x)}^\nu)$.
For $s<t<s+\tau_s(x)$, the transition kernel has an atom at $\psi_{s,t}(x)$ and an absolutely continuous part with the same support as $\mu_t^\nu$:
\begin{equation}\label{eq:outlier_kernel}
\calP^\nu_{s,t}(x,\d y)
=\frac{t-s}{|F^\nu_{s,t}(y)-x|^2}\mu_t^\nu(\d y)
+\pp{1-\frac{t-s}{\tau_s(x)}}\delta_{\psi_{s,t}(x)}(\d y),
\end{equation}
where
\[
F^\nu_{s,t}(y)=\lim_{\eta\downarrow0}F^\nu_{s,t}(y+\i\eta),
\qquad y\in\R.
\]
Both expressions follow by inverting the Cauchy transform. Fix $s<t<s+\tau_s(x)$ and set $a=\psi_{s,t}(x)$. By \citet[Lemma~4 and Proposition~2]{biane97free}, $\psi_{s,t}$ is a local inverse of $F_{s,t}^\nu$ near $x$ (i.e., $F_{s,t}^\nu(\psi_{s,t}(u)) = u$ for $u$ in a small neighborhood of $x$). Consequently,
\[
F^\nu_{s,t}(a)=x
\qmand
(F^\nu_{s,t})'(a)
=\frac1{1+(t-s)G_{\mu_s^\nu}'(x)}
=\frac1{1-(t-s)/\tau_s(x)}.
\]
Using the identity $\lim_{\eta\downarrow0}\i\eta G_\mu(a+\i\eta)=\mu(\{a\})$ and \eqref{eq:resolvent F}, we obtain
\[
\calP^\nu_{s,t}(x,\{a\})
=\lim_{\eta\downarrow0}
\frac{\i\eta}{F^\nu_{s,t}(a+\i\eta)-x}=\lim_{\eta\downarrow0}\frac{\i\eta}{F^\nu_{s,t}(a+\i\eta)-F^\nu_{s,t}(a)}=\frac1{(F^\nu_{s,t})'(a)}=1-\frac{t-s}{\tau_s(x)}.
\]
This proves the probability mass formula in \eqref{eq:outlier_kernel}.

The boundary homeomorphism in \citet[Lemma~4]{biane97free} also shows that $F^\nu_{s,t}(y)=x$ has the unique real solution $y=a$. Since $F^\nu_{s,t}$ extends continuously to the closed upper half-plane, the Cauchy transform $1/(F^\nu_{s,t}(y+\i\eta)-x)$ is locally uniformly bounded for $y$ bounded away from $a$ and all sufficiently small $\eta>0$. Stieltjes inversion therefore implies that $\calP^\nu_{s,t}(x,\cdot)$ is absolutely continuous on $\R\setminus\{a\}$. Thus its singular part consists only of the atom computed in \eqref{eq:outlier_kernel}.

For the absolutely continuous part, recall that its density function $\rho$ can be recovered from the Cauchy transform by $\rho(y)=-\pi\inv\lim_{\eta\downarrow0}\rmIm G_\mu(y+\i\eta)$, for Lebesgue-almost every $y\in\R$. Let $p^\nu_{s,t}(x,y)$ denote the density of the absolutely continuous part of the transition kernel $\calP^\nu_{s,t}(x,\d y)$, and $p^\nu_t(y)$ the density of $\mu_t^\nu$. Using \eqref{eq:resolvent F} and \eqref{eq:F}, we obtain, for Lebesgue-a.e.\ $y\in\R$,
\begin{align*}
p^\nu_{s,t}(x,y)
&=-\frac1\pi\rmIm\frac1{F^\nu_{s,t}(y)-x}
=\frac{\rmIm F^\nu_{s,t}(y)}
{\pi|F^\nu_{s,t}(y)-x|^2},\\
p^\nu_t(y)
&=-\frac1\pi\rmIm\frac{y-F^\nu_{s,t}(y)}{t-s}
=\frac{\rmIm F^\nu_{s,t}(y)}{\pi(t-s)}.
\end{align*}
Combining these identities proves the density formula part in \eqref{eq:outlier_kernel}.
\end{remark}
In particular, we point out that the classical free Meixner processes can all be represented using the standard classical free Brownian transition kernels, with a time shift and a spatial translation; see Appendix~\ref{sec:meixner}.

\section{A martingale problem}\label{sec:martingale}
We let $X^{\nu,\lambda}$ denote a Markov process with transition kernels $(\calP^\nu_{s,t})_{0\le s<t}$ defined in Section \ref{sec:classical free Bm} and initial law $\lambda$, and assume that its paths are in $D([0,\infty))$ (recall Lemma \ref{lem:D}). The main contribution of this section is to show that the process can be characterized as the unique solution of a martingale problem in Theorem \ref{thm:martingale_problem}.
Set
\[
\calF_t^{\nu,\lambda} := \sigma\left(X^{\nu,\lambda}(s):0\leq s\leq t\right), \qquad t\geq0.
\]
Introduce
\equh\label{eq:calW}
\calW := \ccbb{f\in C^2(\R)\cap C_0(\R): f''\in C_0(\R)},
\qquad \nn f_{\calW} := \nn f_\infty+\nn{f''}_\infty.
\eque
\begin{theorem}\label{thm:martingale_problem}
The process $X^{\nu,\lambda}$ is the only solution to the $(\calA_t^\nu,\calW)$-martingale problem  with initial law $\lambda$ in the following sense. 
\begin{enumerate}[(i)]
\item (Existence) For every $f\in\calW$,
\[
M_f^{\nu,\lambda}(t)
:=
f\pp{X^{\nu,\lambda}(t)}-f\pp{X^{\nu,\lambda}(0)}
-\int_0^t \calA_s^\nu f\pp{X^{\nu,\lambda}(s)}\d s,
\qquad t\geq0,
\]
is a martingale with respect to
$(\calF_t^{\nu,\lambda})_{t\geq0}$.

\item (Uniqueness) Let $Y=(Y(t))_{t\geq0}$ be any càdlàg process adapted to
a filtration $(\calG_t)_{t\geq0}$ such that
$\calL(Y(0))=\lambda$. If, for every $f\in\calW$,
\[
f(Y(t))-f(Y(0))
-\int_0^t\calA_s^\nu f(Y(s))\,\d s,
\qquad t\geq0,
\]
is a martingale with respect to $(\calG_t)_{t\geq0}$, then
the laws of $Y$ and $X^{\nu,\lambda}$ are the same as probability measures on $D([0,\infty))$.
\end{enumerate}
\end{theorem}
\begin{proof}
The proof of existence is the same as that of \citet[Lemma~2]{anshelevich13generators}; we include the argument for completeness.
For $0\leq s\leq t$, we have
\begin{align}
\mathbb E&\pp{
M^{\nu,\lambda}_f(t)-M^{\nu,\lambda}_f(s)\mmid\calF^{\nu,\lambda}_s\nonumber
}
\\
& = \esp\pp{f\pp{X^{\nu,\lambda}(t)}\mmid\calF^{\nu,\lambda}_s}- f\pp{X^{\nu,\lambda}(s)} - \esp\pp{\int_s^t\calA^\nu_rf\pp{X^{\nu,\lambda}(r)}\d r\mmid\calF^{\nu,\lambda}_s}\nonumber
\\&=
\calP_{s,t}^\nu f\left(X^{\nu,\lambda}(s)\right)-f\left(X^{\nu,\lambda}(s)\right) 
-\int_s^t\calP_{s,r}^\nu\calA_r^\nu f \left(X^{\nu,\lambda}(s)\right)\d r.\label{eq:=0}
\end{align}
By \citet[Proposition~1]{anshelevich13generators}, for every
$f\in\calW$ and $0\leq s\leq t$,
\[
\calP_{s,t}^\nu f =f+\int_s^t \calP_{s,r}^\nu\calA_r^\nu f\,\d r,
\]
where the integral is understood in $C_0(\R)$. Therefore, the expression in \eqref{eq:=0} is zero and hence $(M^{\nu,\lambda}_f)_{t\geq0}$ is a martingale with respect to $(\calF^{\nu,\lambda}_t)_{t\geq0}$.

We next prove the uniqueness. Let $Y$ and $\calG$ be as in the statement. Fix $t>0$. Let $F_{s,t}^\nu(z)$ be as in \eqref{eq:F}, and for each $z\in\C^+$ set
\[
v(r,x)\equiv v\topp {t,z}(r,x) :=\frac1{F_{r,t}^\nu(z) - x}, \mfa r\in[0,t]. 
\]
The first step is to show that for every $t,z$ fixed, $(v(r,Y_r))_{r\in[0,t]}$ is a martingale.
The function $v$ is complex-valued, and both its real and imaginary parts will be shown below to be elements in $C^1([0,t];\calW)$, the space of continuous $\calW$-valued functions with continuous derivative, equipped with the norm 
\[
\nn u_{C^1([0,t];\calW)} :=\sup_{s\in[0,t]}\nn {u(s,\cdot)}_\calW + \sup_{s\in[0,t]}\nn{\partial_r u(s,\cdot)}_\calW.
\] 
We show that $(v(r,Y_r))_{r\in[0,t]}$ is a martingale in two steps. \medskip

\noindent (i) First, we show that   for all $u\in C^1([0,t];\calW)$, 
\equh\label{eq:martingale}
M^*_u(s):=u(s,Y_s)-u(0,Y_0) - \int_0^s\pp{\partial_ru(r,Y_r)+[\calA_r^\nu u(r,\cdot)](Y_r)}\d r, s\in[0,t],
\eque
is a martingale with respect to $(\calG_s)_{s\in[0,t]}$. (Note that here $\calA_t^\nu$ is a generator, and it is not known that it is the generator of $Y$ yet.) 
This fact is established by an approximation argument. Recall that for $f\in\calW$, by assumption
\[
\wt M_f(s) := f(Y_s)-f(Y_0)-\int_0^s \calA_r^\nu f(Y_r)\d r, s\ge 0,
\]
is a martingale with respect to $\calG$.  
We start with $u(r,x) = a(r)f(x)$, $a\in C^1([0,t]), f\in\calW$, and write
\[
M_u^*(s)=a(s)f(Y_s) - a(0)f(Y_0) - \int_0^s\pp{a'(r)f(Y_r)+a(r)\calA_r^\nu f(Y_r)}\d r, s\in[0,t].
\]
Using
\[
f(Y_s) = f(Y_0)+\wt M_f(s)+\int_0^s\calA_r^\nu f(Y_r)\d r =: f(Y_0)+\wt M_f(s)+V(s),
\]
we have
\begin{align*}
M_u^*(s)& = 
a(s)\pp{f(Y_0)+\wt M_f(s)+V(s)} - a(0)f(Y_0)\\
& \quad - \int_0^s\pp{a'(r)\pp{f(Y_0)+\wt M_f(r)+V(r)}+a(r)\calA_r^\nu f(Y_r)}\d r\\
& = a(s)\wt M_f(s)+a(s)V(s)-\int_0^sa'(r)\pp{\wt M_f(r)+V(r)}+a(r)\calA_r^\nu f (Y_r)\d r\\
&= a(s)\wt M_f(s)-\int_0^sa'(r)\wt M_f(r)\d r,
\end{align*}
where in the last step we used the fact that $V$ is a function of finite variation and applied integration by parts to $a(r)V(r)$. Then, 
one readily checks that $\esp(M_u^*(s)\mid\calG_r) = M_u^*(r)$, and hence $M_u^* = (M_u^*(s))_{s\in[0,t]}$ is a martingale with $u(s,x) = a(s)f(x)$. (Or, by integration by parts for semimartingales one can write directly $M_u^*(s) =  \int_0^s a(r)\d \wt M_f(r)$.)

Now we check that for every $u\in C^1([0,t];\calW)$, $M_u^*$ is a martingale. We first consider
\[
u_n(r,\cdot) :=\summ k0n\binom nk\pp{\frac rt}^k\pp{1-\frac rt}^{n-k}u\pp{\frac{kt}n,\cdot}, r\in[0,t], n\in\N,
\]
and they satisfy $u_n\to u$ in $C^1([0,t];\calW)$ as $n\to\infty$ (this is the Banach-valued Bernstein approximation). For each $n\in\N$, $u_n$ is of the form $\summ i1d a_i(r)f_i(x)$, and hence $M_{u_n}^*$ is a martingale.
We next show that 
\[
\limn\nn{M^*_{u_n}- M^*_u}_\infty= 0,
\] uniformly for all $\omega\in\Omega$, with $\nn g_\infty :=\sup_{s\in[0,t]}\abs{g(s)}$. Indeed, 
\begin{multline*}
\abs{M^*_{u_n}(s)-M^*_u(s)} \le \abs{u_n(s,Y_s)-u(s,Y_s)}+\abs{u_n(0,Y_0)-u(0,Y_0)} \\
 +\int_0^s \abs{\partial_ru_n(r,Y_r)-\partial_ru(r,Y_r)}+\abs{\bb{\calA_r^\nu u_n(r,\cdot)}(Y_r)-\bb{\calA_r^\nu u(r,\cdot)}(Y_r)}\d r.
\end{multline*}
Since $\nn{\calA_r^\nu g}_{\infty}\le  \nn g_\calW/2$, 
we have
\begin{align*}
\nn{M^*_{u_n}-M^*_u}_\infty
& \le \pp{2+t+\frac t2}\nn{u_n-u}_{C^1([0,t];\calW)}.
\end{align*}
In particular, it follows that
\[
\mathbb E\pp{M^*_u(s)\mmid\calG_r}
=\lim_{n\to\infty}\mathbb E\pp{M^*_{u_n}(s)\mmid\calG_r}
=\lim_{n\to\infty}M^*_{u_n}(r)=M^*_u(r), \mbox{ almost surely,}
\]
and hence $M^*_u$ is a martingale for all $u\in C^1([0,t];\calW)$.

Set $\calW_\C:=\{f: \mathrm {Re}f\in \calW, \mathrm{Im}f\in\calW\}$ and $\nn f_{\calW_\C}:=\nn{\mathrm{Re}f}_\calW+\nn{\mathrm{Im}f}_\calW$. So the above development yields that for every $u\in C^1([0,t];\calW_\C)$, $M_u^*$ is a complex-valued martingale. 
 We now check $v\in C^1([0,t];\calW_\C)$. Recall 
\[
v(r,x) = \frac1{F_{r,t}^\nu(z)-x} = \frac1{z-(t-r)G_{\mu^\nu_t}(z)-x}=: \frac1{w(r)-x},
\] 
and that $t>0$ and $z\in\C^+$ are fixed. Observe that  $z\in \C^+$  and $G_{\mu^\nu_t}(z)\in \C^-$, and hence the imaginary part of the denominator of $v(r,x)$, $F_{r,t}^\nu(z)-x$,  is uniformly bounded away from zero for all $x\in\R, r\in[0,t]$. One readily checks that $v(r,\cdot)$ is in $\calW_\C$, $\partial_rv(r,\cdot) = - w'(r)/(w(r)-x)^2\in\calW_\C$, and $r\mapsto\partial_r v(r,\cdot)$ is continuous with respect to $\nn\cdot_{\calW_\C}$. So $v\in C^1([0,t];\calW_\C)$ and hence $M^*_v$ is a martingale for all $t>0$. 
\medskip

\noindent (ii) Second, we show that 
\equh\label{eq:partial_r+A_r}
\partial_r v(r,x) + [\calA_r^\nu v(r,\cdot)](x) = 0.
\eque
Indeed, with $h_z(x):=1/(z-x)$, we have
\[
\partial_x\pp{\frac{h_z(x)-h_z(y)}{x-y}} = \frac1{z-y}\frac1{(z-x)^2}.
\]
Therefore
\[
\calA_r^\nu h_z(x) = \frac1{(z-x)^2}\int\frac{\mu_r^\nu(\d y)}{z-y} = \frac{G_{\mu_r^\nu}(z)}{(z-x)^2},
\]
and 
\[
\bb{\calA_r^\nu v(r,\cdot)}(x) = \frac{G_{\mu_r^\nu}(F_{r,t}^\nu(z))}{(F_{r,t}^\nu(z)-x)^2} = \frac{G_{\mu_t^\nu}(z)}{(F_{r,t}^\nu(z)-x)^2}.
\]
Moreover, 
\[
\partial_r v(r,x) = \partial_r\pp{\frac1{z-(t-r)G_{\mu_t^\nu}(z)-x}} = -\frac{G_{\mu_t^\nu}(z)}{(F_{r,t}^\nu(z)-x)^2}.
\]
So \eqref{eq:partial_r+A_r} follows. 

The identity \eqref{eq:partial_r+A_r}  and the claim that $M_v^*(s)$ in \eqref{eq:martingale} is a martingale together lead to the conclusion that $(v(s,Y_s)-v(0,Y_0))_{s\in[0,t]}$ is a martingale, and hence so is $(v(s,Y_s))_{s\in[0,t]}$. This implies in particular
\equh\label{eq:Y resolvent}
\esp\pp{\frac1{z-Y_t}\mmid\calG_s} = \esp (v(t,Y_t)\mid\calG_s) = v(s,Y_s) = \frac1{F_{s,t}^\nu(z)-Y_s}, \mbox{ almost surely,}
\eque
for all $0\le s<t$. \medskip

With \eqref{eq:Y resolvent} established, we can conclude that $Y$ is Markov with transition kernels $\calP^\nu$. Indeed, setting $h_z(x):=(z-x)^{-1}$, \eqref{eq:Y resolvent} yields, for $0\leq s<t$ and $z\in\mathbb C^+$,
\[
\mathbb E\bigl(h_z(Y_t)\mid\calG_s\bigr)
=
\calP_{s,t}^\nu h_z(Y_s), \mbox{ almost surely.}
\]
By complex conjugation, the same identity holds for $z\in\mathbb C^-$. The complex linear span $\calA := \mathrm{span}\{h_z:z\in\mathbb C\setminus\mathbb R\}$ is dense in $C_0(\mathbb R;\mathbb C)$ with respect to the uniform norm $\nn\cdot_\infty$ by the Stone--Weierstrass theorem. Indeed, $\wb h_z=h_{\wb z}$, so $\calA$, and hence its uniform closure $\wb\calA$, is closed under complex conjugation. Moreover, $\calA$ separates points and vanishes nowhere. Its uniform closure $\wb {\calA}$ is also an algebra. To see this, we first notice that for $z\ne w$, $h_z(x)h_w(x) = (h_z(x)-h_w(x))/(w-z)$ and $h_z^2$ can be approximated by $h_zh_{w_n}$ with $w_n\to z$ with respect to the uniform norm. So $h_zh_w\in\wb\calA$ for all $z,w\in\C\setminus\R$, and hence for all $f,g\in\calA$, $fg\in\wb\calA$. One can then check directly that $\wb\calA$ is an algebra. The closure $\wb\calA$ therefore satisfies the assumptions of the Stone--Weierstrass theorem \citep[Theorems~4.51--4.52]{folland99real}, which implies that $\wb\calA=C_0(\mathbb R;\mathbb C)$.

Since conditional expectation and $\calP_{s,t}^\nu$ are contractions,
it follows that
\[
\mathbb E\bigl(g(Y_t)\mid\calG_s\bigr)
=
\calP_{s,t}^\nu g(Y_s), \mbox{ almost surely, for every }
 g\in C_0(\mathbb R;\C).
\]
A monotone class argument then extends this identity to every bounded Borel function $g$. Thus $Y$ is a time-inhomogeneous Markov process,
with respect to $(\calG_t)$, with transition kernels $(\calP_{s,t}^\nu)_{0\le s< t}$. It then follows that $Y$ and $X^{\nu,\lambda}$ have the same finite-dimensional laws. Both processes have realizations in $D([0,\infty))$, and hence both have the same law on $D([0,\infty))$. 
\end{proof}


\section{Proofs of the limit theorems}\label{sec:proofs}
We prove Theorem \ref{thm:2} in Section \ref{sec:proof1}. We deal with the Hermite root flow in Section \ref{sec:hermite entrance}. 
\subsection{Proof of Theorem \ref{thm:2}}\label{sec:proof1}
We first compute the generator of $X_n(t) = r_{n,I_n(t)}(t)$. 
\begin{lemma}\label{lem:generator}
The generator of $X_n(t) = r_{n,I_n(t)}(t)$ is, for all $f\in\calW$, 
\equh\label{eq:calAf 2}
\calA_{n,t}f(r_{n,i}(t)) = \int K_f(r_{n,i}(t),y)\mu_{n,t}(\d y)  - \frac1{2n}f''(r_{n,i}(t)), t>0, i=1,\dots,n.
\eque
\end{lemma}
\begin{proof}
We start by writing
\begin{align*}
\esp& \pp{f(X_n(t+h))-f(X_n(t))\mmid X_n (t) = r_{n,i}(t)} \\
&\quad\quad= \esp\pp{f(r_{n,I_n(t+h)}(t+h))-f(r_{n,I_n(t+h)}(t))\mmid I_n (t) = i}\\
& \quad\quad\quad + \esp \pp{f(r_{n,I_n(t+h)}(t))-f(r_{n,I_n(t)}(t))\mmid I_n(t) = i}.
\end{align*}
First, 
\begin{align*}
\lim_{h\downarrow0}\frac1h & \esp \pp{f(r_{n,I_n(t+h)}(t+h))-f(r_{n,I_n(t+h)}(t))\mmid I_n(t) = i}\\
& = \lim_{h\downarrow0}\frac1h\esp \pp{f(r_{n,I_n(t+h)}(t+h))-f(r_{n,I_n(t+h)}(t))\mmid I_n(t+h) = i} \\
& = f'(r_{n,i}(t))r_{n,i}'(t) = -\frac1n\sum_{j\ne i}\frac{f'(r_{n,i}(t))(r_{n,j}(t)-r_{n,i}(t))}{(r_{n,j}(t)-r_{n,i}(t))^2},
\end{align*}
where in the last step we recalled the ODE for $r_{n,i}(t)$ in  \eqref{eq:frozen DBm}.
Second, 
recalling the generator of the process $I_n$ in \eqref{eq:I generator}, we have
\begin{align*}
\lim_{h\downarrow0}\frac1h\esp \pp{f(r_{n,I_n(t+h)}(t))-f(r_{n,I_n(t)}(t))\mmid I_n (t) = i}& = \bb{\calL\topp n_t f(r_{n,\cdot}(t))}(i)\\
& = \frac1n\sum_{j\ne i}\frac{f(r_{n,j}(t))-f(r_{n,i}(t))}{(r_{n,j}(t)-r_{n,i}(t))^2}.
\end{align*}
Introduce \equh\label{eq:K_f}
K_f(x,y) := \begin{cases}
\displaystyle \frac{f(y)-f(x)-f'(x)(y-x)}{(y-x)^2}, &\mbox{ if } y\ne x,\\\\
\displaystyle \frac12 f''(x), &\mbox { if } x=y,
\end{cases}
\eque
where the case $x=y$ is defined by continuity. We have shown that 
\begin{align}
\calA_{n,t}f(r_{n,i}(t))& = \lim_{h\downarrow 0}\frac1h\esp \pp{f(X_n(t+h))-f(X_n(t))\mid X_n(t) = r_{n,i}(t)} \nonumber\\
&= \frac1n \sum_{j\ne i}K_f(r_{n,i}(t),r_{n,j}(t)).
\label{eq:calAf 1}
\end{align}
(The convergence is also uniform over $i=1,\dots,n$.) The claimed formula follows by recalling that  $\mu_{n,t} = n\inv\summ i1n \delta_{r_{n,i}(t)}$.
\end{proof}
Next, we shall set up the problem as a martingale problem. Recall $\calW$ defined in \eqref{eq:calW}.
We start by focusing on the process $X_n$. Consider $\calF_t\topp n := \sigma(I_n(s): s\in[0,t])$. 
\begin{lemma}\label{lem:M_n^f}
For $f\in\calW$, the process 
\[
M_n^f(t):=f(X_n(t))-f(X_n(0)) - \int_0^t \calA_{n,s}f(X_n(s-))\d s, t\ge 0,
\]
is a martingale with respect to the filtration $(\calF_t\topp n)_{t\ge 0}$. 
\end{lemma}
\begin{proof}Note that $(X_n(t))_{t\ge 0}$ is a piecewise deterministic Markov process. In particular, 
we have for all $f\in\calW$,
\begin{multline*}
f(X_n(t))-f(X_n(0)) \\
= \int_0^t f'(r_{n,I_n(s-)}(s))r'_{n,I_n(s-)}(s)\d s + \sum_{s\in(0,t]: I_n(s-)\ne I_n(s)}\pp{f(X_n(s))-f(X_n(s-))}.
\end{multline*}
For $i\ne j$, let $N_{n,i,j}(t)$ denote the counting process of the jumps of the process $I_n$ from $i$ to $j$. That is,
\[
N_{n,i,j}(t) := \sum_{s\in(0, t]} \inddd{I_n(s-) = i, I_n(s) = j}.
\]
Accordingly
\[
\wb N_{n,i,j}(t) := N_{n,i,j}(t) - \int_0^t \inddd{I_n(s-) = i}Q_{i,j}\topp n(\vec r_n(s))\d s,
\]
is a local martingale with respect to $(\calF_t\topp n)_{t\ge 0}$ \citep[Theorem II.T8, p.~27]{bremaud81point}. Moreover, it is square-integrable as for every $n,i,j$ and $t>0$, we have
\[
\sup_{s\in[0,t]}Q\topp n_{i,j}(\vec r_n(s))<\infty,
\]
because of the assumption that the initial roots are all distinct.  Therefore, 
\begin{align*}
\sum_{s\in(0,t]: I_n(s-)\ne I_n(s)}&\pp{f(X_n(s))-f(X_n(s-))} = \summ i1n \sum_{j\ne i} \int_0^t \pp{f(r_{n,j}(s))-f(r_{n,i}(s))}N_{n,i,j}(\d s)\\
& = \summ i1n\sum_{j\ne i} \bigg(\int_0^t \pp{f(r_{n,j}(s))-f(r_{n,i}(s))}\wb N_{n,i,j}(\d s) \\
& \quad+ \int_0^t \pp{f(r_{n,j}(s))-f(r_{n,i}(s))}\inddd{I_n(s-) = i}Q_{i,j}\topp n(\vec r_n(s))\d s\bigg)\\
& = \summ i1n\sum_{j\ne i} \int_0^t \pp{f(r_{n,j}(s))-f(r_{n,i}(s))}\wb N_{n,i,j}(\d s) \\
&\quad+  \int_0^t\sum_{j\ne I_n(s-)} \pp{f(r_{n,j}(s))-f(X_n(s-))}Q_{I_n(s-),j}\topp n(\vec r_n(s))\d s.
\end{align*}
In the last expression above, the first double sum is a square-integrable martingale (in $t$) since $|f(r_{n,j}(s))-f(r_{n,i}(s))|\le 2\nn f_\infty$. The second term can be rewritten as (recalling \eqref{eq:calAf 1})
\[
\int_0^t \calA_{n,s}f(X_n(s-))\d s - \int_0^t f'(r_{n,I_n(s-)}(s))r_{n,I_n(s-)}'(s)\d s.
\]
Combining all the above, we have shown that
\[
f(X_n(t))-f(X_n(0)) = \int_0^t\calA_{n,s}f(X_n(s-))\d s + \summ i1n \sum_{j\ne i}\int_0^t(f(r_{n,j}(s))-f(r_{n,i}(s)))\wb N_{n,i,j}(\d s).
\]This completes the proof. 
\end{proof}

Next, we collect a couple of estimates. 

\begin{lemma}\label{lem:generator_convergence}
Suppose that $n\inv\summ i1n\delta_{r_{n,i}(0)}\weakto\nu$. Then, for every $f\in\calW$ and $t>0$,
\[
\limn\max_{1\le i\le n}
\left|\calA_{n,t}f(r_{n,i}(t))-\calA_t^\nu f(r_{n,i}(t))\right|=0.
\]
Moreover, for every $T>0$,
\[
\limn\int_0^T\max_{1\le i\le n}
\left|\calA_{n,t}f(r_{n,i}(t))-\calA_t^\nu f(r_{n,i}(t))\right|\d t=0.
\]
\end{lemma}

\begin{proof}
By \eqref{eq:calAf 2} and the identity
\[
\calA_t^\nu f(x)=\int_\R K_f(x,y)\mu_t^\nu(\d y),
\]
we have
\begin{multline*}
\calA_{n,t}f(r_{n,i}(t))-\calA_t^\nu f(r_{n,i}(t)) \\= \int_\R K_f(r_{n,i}(t),y)\mu_{n,t}(\d y) - \int_\R K_f(r_{n,i}(t),y)\mu_t^\nu(\d y) - \frac1{2n}f''(r_{n,i}(t)).
\end{multline*}
Since $f''\in C_0(\R)$, the last term goes to $0$ as $n\to\infty$ uniformly for $i=1,\dots,n$. It remains to control the difference of the first two integrals.

Taylor's formula gives
\[
K_f(x,y)=\int_0^1(1-u)f''((1-u)x+uy)\d u,
\qquad
\|K_f\|_\infty\le\frac12\|f''\|_\infty.
\]
We first show that
\equh\label{eq:AA}
\limn\sup_{x\in\R}
\left|\int_\R K_f(x,y)\mu_{n,t}(\d y) -\int_\R K_f(x,y)\mu_t^\nu(\d y)\right|
=0.
\eque
To see this, for every $\epsilon>0$, first take $B_R = [-R,R]$ with $R>0$ such that $\sup_{n\in\N}\mu_{n,t}(B_R^c)+\sup_{n\in\N}\mu_t^\nu(B_R^c)<\epsilon$. Next, since $f''\in C_0(\R)$, it is uniformly continuous, and consequently, the family $\{K_f(x,\cdot)\}_{x\in\R}$ is uniformly bounded and uniformly equicontinuous. Moreover, the family of functions $\{K_f(x,\cdot)\}_{x\in\R}$ restricted to $B_R$ is relatively compact in $C(B_R)$ with respect to the uniform norm by the Arzel\`a--Ascoli theorem. In particular, there exist $m$ (depending on $\epsilon$) and $x_1,\dots,x_m$ such that for all $x\in\R$, there exists $i=1,\dots,m$ such that
\[
\sup_{y\in B_R}\abs{K_f(x,y)-K_f(x_i,y)}<\epsilon. 
\]
For $x$ and $x_i$ as above, we have that
\begin{multline*}
\abs{\int_{B_R}K_f(x,y)\mu_{n,t}(\d y) - \int_{B_R}K_f(x,y)\mu_t^\nu(\d y)}\\
 \le 
\abs{\int_{B_R}K_f(x_i,y)\mu_{n,t}(\d y) - \int_{B_R}K_f(x_i,y)\mu_t^\nu(\d y)}  \\
+ \int_{B_R}\abs{K_f(x,y)-K_f(x_i,y)}(\mu_{n,t}(\d y)+\mu_t^\nu(\d y)).
\end{multline*}
The first term on the right-hand side above goes to zero because of $\mu_{n,t}\weakto\mu_t^\nu$ as $n\to\infty$. The second term is bounded by $2\epsilon$. Combining all the above, we have shown that
\[
\limsupn \sup_{x\in\R}
\left|\int_\R K_f(x,y)\mu_{n,t}(\d y) -\int_\R K_f(x,y)\mu_t^\nu(\d y)\right|\le (2+\nn{f''}_\infty)\epsilon. 
\]
Since $\epsilon>0$ is arbitrary, we have proved \eqref{eq:AA}, and the first assertion. Finally,
\[
\max_{1\le i\le n}|\calA_{n,t}f(r_{n,i}(t))|
\le\frac12\|f''\|_\infty,
\qquad
\|\calA_t^\nu f\|_\infty\le\frac12\|f''\|_\infty.
\]
The second assertion follows from the dominated convergence theorem.
\end{proof}
Next, we prove the tightness of $(X_n)_{n\in\N}$. 
\begin{lemma}\label{lem:tightness}
Suppose that the initial laws $(\calL(X_n(0)))_{n\in\N}$ form a tight family. Then the laws of $(X_n)_{n\in\N}$ are tight in $D([0,\infty))$ equipped with the Skorokhod $J_1$ topology.
\end{lemma}

\begin{proof}
For each fixed $n$ and $S>0$, all root trajectories are bounded on $[0,S]$. We may therefore apply Lemma~\ref{lem:M_n^f} to smooth compactly supported functions agreeing with $x$ and $x^2$, respectively, on an interval containing these trajectories. Direct calculation gives
\[
\bb{\calA_{n,t} x}(\cdot)=0,
\qquad
\bb{\calA_{n,t}x^2}(\cdot)=1-\frac1n.
\]
It follows that $(X_n(t))_{t\ge0}$ and $(X_n(t)^2-(1-n\inv)t)_{t\ge0}$ are martingales. In particular, $(X_n(t)-X_n(0))_{t\ge0}$ is a square-integrable martingale with
\[
\left\langle X_n-X_n(0)\right\rangle_t
=\left(1-\frac1n\right)t.
\]
Fix $T>0$. By Doob's maximal inequality, for every $x>0$,
\[
\begin{split}
\proba\left(\sup_{0\le t\le T}|X_n(t)|>x\right)
&\le \proba\pp{|X_n(0)|>\frac x2}
+\proba\left(\sup_{0\le t\le T}|X_n(t)-X_n(0)|>\frac x2\right)\\
&\le \proba\pp{|X_n(0)|>\frac x2}+\frac{4T}{x^2}.
\end{split}
\]
Thus, since $(X_n(0))_{n\in\N}$ is tight, it follows that
\equh\label{eq:Aldous1}
\lim_{x\to\infty}\limsupn \proba\pp{\sup_{t\in[0,T]}|X_n(t)|>x} = 0.
\eque
Moreover, for every $(\calF_t\topp n)$-stopping time $\tau_n\le T$ and every $0\le\tau\le1$, optional sampling gives
\[
\esp\left|X_n(\tau_n+\tau)-X_n(\tau_n)\right|^2
=\left(1-\frac1n\right)\tau\le\tau,
\]
and hence for all $\eta>0$, 
\equh\label{eq:Aldous2}
\lim_{\tau\downarrow 0}\limsupn\proba\pp{|X_n(\tau_n+\tau) - X_n(\tau_n)|>\eta} = 0.
\eque
With \eqref{eq:Aldous1} and \eqref{eq:Aldous2} established, the tightness of the sequence now follows from Aldous' condition \citep[Theorem 16.10]{billingsley99convergence}.
\end{proof}
\begin{proof}[Proof of Theorem \ref{thm:2}]
(i) We first prove that
\equh\label{eq:marginal conv}
X_n\weakto X^{\nu,\lambda}
\eque
as $n\to\infty$ in  $D([0,\infty))$ equipped with the Skorokhod $J_1$ topology.

Since $X_n(0)\weakto\lambda$, Lemma~\ref{lem:tightness} implies
tightness of $(X_n)_{n\in\N}$. So it suffices to show that every subsequential limit has the same law as $X^{\nu,\lambda}$.  Let $Y$ be a subsequential limit in $D([0,\infty))$.
Relabeling this subsequence and applying the Skorokhod representation theorem, we may assume that
\[
X_n \to Y
\quad\text{almost surely in }D([0,\infty)),
\]
as $n\to\infty$. To show that $Y$ has the law of $X^{\nu,\lambda}$, it suffices to show that $Y$ solves the $(\calA_t^\nu,\calW)$-martingale problem as explained in Section \ref{sec:martingale}. Indeed, since $X^{\nu,\lambda}$ has initial law $\lambda$ and solves the same martingale problem, of which the solution is unique following Theorem \ref{thm:martingale_problem}, the law of $Y$ has to be the same as the law of $X^{\nu,\lambda}$. 

We next  show that $Y$ solves the $(\calA_t^\nu,\calW)$-martingale problem.  That is for every $f\in\calW$, 
\equh\label{eq:M^f}
M^f(t):=f(Y(t))-f(Y(0)) -\int_0^t\calA_s^\nu f(Y(s))\d s, \quad t\ge 0,
\eque
is a martingale with respect to the natural filtration of $(Y(t))_{t\ge 0}$. Fix $f\in\calW$. Lemma \ref{lem:M_n^f} yields that 
\[
M_n^f(t)=f(X_n(t))-f(X_n(0))-\int_0^t\calA_{n,s}f(X_n(s))\d s,\quad t\ge 0,
\]
is a martingale with respect to the natural filtration of $(X_n(t))_{t\ge0}$. Note that we write $X_n(s)$ instead of $X_n(s-)$; this does not change the value of the  integral on the right-hand side above. By Lemma~\ref{lem:generator_convergence}, for every $T>0$,
\[
\limn \int_0^T \left|\calA_{n,s}f(X_n(s))-\calA_s^\nu f(X_n(s))\right|\d s = 0.
\]
Since $Y$ is in $D([0,\infty))$, we have $\limn X_n(t) = Y(t)$ for Lebesgue almost every $t\in[0,\infty)$ (in fact, except at countably many points at which $Y$ is discontinuous).  Since $x\mapsto\calA_s^\nu f(x)$ is continuous and
\[
\|\calA_s^\nu f\|_\infty\le\frac12\|f''\|_\infty,
\]
the dominated convergence theorem yields 
\[
\limn\int_0^T|\calA_s^\nu f(X_n(s)) - \calA_s^\nu f(Y(s))|\d s = 0, \quad \mbox{ almost surely,}
\]
and hence 
\[
\limn \sup_{0\le t\le T}
\abs{\int_0^t\calA_{n,s}f(X_n(s))\d s -\int_0^t\calA_s^\nu f(Y(s))\d s }
=0,
\quad\text{almost surely}.
\]
Introduce 
\[
J:=\{t\in(0,\infty):\proba(Y(t-)\ne Y(t))>0\},
\]
which is at most countable. 
Pick a countable set $\mathbb T\subset[0,\infty)\setminus J$ such that $0\in \mathbb T$ and it is dense in $[0,\infty)$. Since $\limn X_n(t) = Y(t)$ for all $t\in\mathbb T$ almost surely, we have thus proved that 
\[
\limn M_n^f(t) = M^f(t), \mfa t\in\mathbb T,
\]
almost surely.

Next, for $t_1,\dots,t_m,s,t\in \mathbb T$ with $0\le t_1<\cdots<t_m\le s<t$ and $H\in C_b(\R^m)$, the martingale
property gives
\[
\esp\pp{\pp{M_n^f(t)-M_n^f(s)}H(X_n(t_1),\dots,X_n(t_m))}=0.
\]
Moreover, for every $T>0$,
\[
\sup_{0\le u\le T}\abs{M_n^f(u)}\le 2\|f\|_\infty+\frac T2\|f''\|_\infty.
\]
The dominated convergence theorem implies that

\[
\esp\pp{\pp{M^f(t)-M^f(s)}
H(Y(t_1),\dots,Y(t_m))}=0.
\]
Fix $s<t, s,t\in\mathbb T$, and consider $\calG_s = \sigma(Y_r: r\le s) = \sigma(Y_r:r\in\mathbb T\cap[0,s])$ (the second equality follows from the right continuity of the paths of $Y$). By a monotone class argument \citep[Theorem 6.3]{jacob03probability}, for every bounded $\calG_s$-measurable random variable $Z$, we have $\esp((M^f(t)-M^f(s))Z) = 0$. This implies that $\esp(M^f(t)\mid\calG_s) = M^f(s)$ almost surely, and hence the identity holds almost surely for all $s<t$ in $\mathbb T$. By right continuity of the paths of $Y$ again, and the fact that $\sup_{t\in[0,T]}|M^f(t)|\le 2\nn f_\infty+(T/2)\nn {f''}_\infty$, the dominated convergence theorem implies
\[
\esp\pp{M^f(t)\mmid\sigma(Y(r):r\in[0,s])} = M^f(s), \mbox{ for all $0\le s<t<\infty$, almost surely.}
\]
That is $(M^f(t))_{t\ge 0}$ defined in \eqref{eq:M^f} is a martingale with respect to $(\calG_t)_{t\ge0}$. Since $f\in\calW$ was arbitrary, $Y$ solves the $(\calA_t^\nu,\calW)$-martingale problem with initial law $\lambda$.\bigskip

\noindent (ii) We prove the joint convergence \eqref{eq:joint}. Recall that if $\nu$ has no atoms, we set $\epsilon=0$; otherwise pick any $\epsilon>0$. Now, \eqref{eq:joint} shall follow from
\equh\label{eq:joint'}
\pp{\pp{\Phi_t^\nu(X_n(t))}_{t\ge\epsilon},\pp{X_n(t)}_{t\ge0}}
\weakto
\pp{\pp{\Phi_t^\nu(X^{\nu,\lambda}(t))}_{t\ge\epsilon},\pp{X^{\nu,\lambda}(t)}_{t\ge0}}
\eque
in $D([\epsilon,\infty))\times D([0,\infty))$ as $n\to\infty$, and
\equh\label{eq:joint2}
\limn\sup_{t\in[\epsilon,T]}\abs{\frac{I_n(t)}n-\Phi_t^\nu(X_n(t))} = 0, \mfa T>\epsilon.
\eque

We first prove \eqref{eq:joint'}. We first show that $(t,y)\mapsto\Phi_t^\nu(y)$ is jointly continuous from $[\epsilon,\infty)\times\R$ to $[0,1]$. 
Since $X^{\nu,\nu}$ is stochastically continuous,  we have $X^{\nu,\nu}(t)\weakto X^{\nu,\nu}(s)$ as $t\to s$ for $s\ge \epsilon$. This weak convergence and the fact that $X^{\nu,\nu}(s)$ has no atoms imply that 
\[
\lim_{t\to s} \sup_{y\in\R}\abs{\Phi_{t}^\nu(y)-\Phi_s^\nu(y)} = 0.
\]
This result is sometimes referred to as P\'olya's theorem. In particular, with $t_n,t\in[\epsilon,T], t_n\to t, y_n\to y$, we have
\[
\abs{\Phi_{t_n}^\nu(y_n)-\Phi_t^\nu(y)} \le \abs{\Phi_{t_n}^\nu(y_n)-\Phi_t^\nu(y_n)}+\abs{\Phi_{t}^\nu(y_n)-\Phi_t^\nu(y)}\to 0
\]
as $n\to\infty$.

We next show that 
\equh\label{eq:Phi(x)}
(x(t))_{t\ge0}\mapsto\left(\Phi_t^\nu(x(t))\right)_{t\ge\epsilon}
\eque
is continuous at every $(x(t))_{t\ge0}$ in $D([0,\infty))$ that is continuous at $\epsilon$ for the Skorokhod $J_1$ topology. Indeed, to show this it suffices to show the convergence of the functions restricted to $[\epsilon,T]$ for every $T>\epsilon$ such that $x$ is continuous at $T$. Fix $x$, such a $T$, and pick a sequence $x_n\to x$ in $D([0,\infty))$. Then, there exists a sequence $(\lambda_n^T)_{n\in\N}$ of continuous and strictly increasing bijections from $[\epsilon,T]$ to $[\epsilon,T]$ such that 
\[
\limn\sup_{t\in[\epsilon,T]}\abs{\lambda_n^T(t)-t} + \sup_{t\in[\epsilon,T]}\abs{x_n\pp{\lambda_n^T(t)}-x(t)} = 0.
\]
(See \citet[Theorem 16.2]{billingsley99convergence}.) 
This implies that 
\[
\limn \abs{\Phi_{\lambda_n^T(t)}^\nu\pp{x_n\pp{\lambda_n^T(t)}}-\Phi_t^\nu(x(t))} = 0
\]
for every $t\in[\epsilon,T]$. By uniform continuity of $(t,y)\mapsto \Phi_t^\nu(y)$ over a compact subset of $[\epsilon,\infty)\times\R$, say $[\epsilon,T]\times [-K,K]$ (e.g., taking $K = \sup_{t\in[\epsilon,T]}|x(t)|+1$ we have for $n$ large enough $\sup_{t\in[\epsilon,T]}|x_n(t)|\le K$), it follows that the above convergence holds uniformly for $t\in[\epsilon,T]$, and hence $(\Phi_t^\nu(x_n(t)))_{t\in[\epsilon,T]}$ converges to $(\Phi_t^\nu(x(t)))_{t\in[\epsilon,T]}$ in $D([\epsilon,T])$. 

We have shown that the mapping in \eqref{eq:Phi(x)} is continuous with respect to the $J_1$ topology at every $(x(t))_{t\ge0}$ that is continuous at $t=\epsilon$. Since $X^{\nu,\lambda}$ is stochastically continuous, it is almost surely continuous at $t=\epsilon$, and hence the continuous mapping theorem yields \eqref{eq:joint'}.

We next prove \eqref{eq:joint2}. 
Let $(I_n^*(t))_{t\ge 0}$ denote the Markov process on $\{1,\dots,n\}$ with uniform initial distribution on $\{1,\dots,n\}$. Recall that this is the stationary distribution of the process and hence the marginal law of the process remains unchanged. Set $X_n^*(t) = r_{n,I_n^*(t)}(t)$. Then, 
\[
X_n^*\weakto X^{\nu,\nu} \quad\text{in }D([0,\infty))
\]
by the first part of the proof of \eqref{eq:marginal conv}. 
Set
\[
\Phi_{n,t}^\nu(y):=\proba(X_n^*(t)\le y).
\]
Then, the facts that $X_n^*(t)\weakto X^{\nu,\nu}(t)$ and also that the law of $X^{\nu,\nu}(t)$ has no atoms imply that, by P\'olya's theorem again,
\[
\limn \sup_{y\in\R}\abs{\Phi_{n,t}^\nu(y)-\Phi_t^\nu(y)} = 0, \mfa t>0. 
\]
We shall, however, need the following stronger result:
\begin{equation}\label{eq:uniform conv}
\limn\sup_{t\in[\epsilon,T]}\sup_{x\in\R}\abs{\Phi_{n,t}^\nu(x)-\Phi_{t}^\nu(x)}=0, \mfa T>\epsilon.
\end{equation}
We prove by contradiction. Suppose the above does not hold. Then, there exist sequences $(n_k)_{k\in\N}, (t_k)_{k\in\N}$ and $\eta>0$ such that $n_k\uparrow\infty$ and 
\equh\label{eq:counter}
\sup_{y\in\R}|\Phi_{n_k,t_k}^\nu(y)-\Phi_{t_k}^\nu(y)|>\eta \mfa k\in\N.
\eque By taking a subsequence and relabeling, we assume further that $t_k\in[\epsilon,T], t_k\to t\in[\epsilon,T]$. But we also have
\equh\label{eq:counter 1}
\lim_{k\to\infty}\sup_{y\in\R}\abs{\Phi_{t_k}^\nu(y)-\Phi_t^\nu(y)} = 0
\eque
because of $X^{\nu,\nu}(t_k)\weakto X^{\nu,\nu}(t)$ and P\'olya's theorem, and 
\equh\label{eq:counter 2}
\lim_{k\to\infty}\sup_{y\in\R}\abs{\Phi_{n_k,t_k}^\nu(y)-\Phi_t^\nu(y)} = 0.
\eque
Indeed, to see \eqref{eq:counter 2} we recall that $X^*_{n_k}\weakto X^{\nu,\nu}$ in $D[0,\infty)$ and $X^{\nu,\nu}$ is stochastically continuous; these two facts yield $X^*_{n_k}(t_k)\weakto X^{\nu,\nu}(t)$ from which \eqref{eq:counter 2} follows from P\'olya's theorem. But the above two limit statements \eqref{eq:counter 1} and \eqref{eq:counter 2} with the triangle inequality contradict \eqref{eq:counter}. We have thus proved \eqref{eq:uniform conv}.

Since the roots are distinct and increasingly ordered, a key observation here is that 
\[
\frac{I_n(t)}n=\Phi_{n,t}^\nu(X_n(t)).
\]
Therefore, by \eqref{eq:uniform conv}, we have
\[
\sup_{t\in[\epsilon,T]} \left|\frac{I_n(t)}n-\Phi_t^\nu(X_n(t))\right|\le\sup_{t\in[\epsilon,T]}\sup_{y\in\R}
|\Phi_{n,t}^\nu(y)-\Phi_t^\nu(y)|\to0
\]
for every $T>\epsilon$. We have thus established \eqref{eq:joint2}.
\end{proof}

\begin{remark}\label{rem:zero}
When $\nu$ has atoms, and $\lambda(\{x: \nu(\{x\})>0\})>0$, the limiting rank process admits no c\`adl\`ag extension to $t=0$. Indeed, consider $U_t=\Phi_t^\nu(X^{\nu,\nu}(t))$, which is uniformly distributed over $(0,1)$ for every $t>0$. Pick $a$ such that  $m=\nu(\{a\})>0$ and $\lambda(\{a\})>0$. By conditioning on the initial value, for every interval $B\subset[0,1]$,
\[
\proba\pp{U_t\in B\mmid X^{\nu,\nu}(0)=a} \le \frac{|B|}{m},\qquad t>0.
\]
Under this conditional law, Blumenthal's zero-one law, applied to the associated space--time Feller process $(t,X^{\nu,\nu}_t)_{t\ge0}$ started at $(0,a)$, implies that the event that $\lim_{t\downarrow0}U_t$ exists has probability zero or one. If it had probability one, the limit would be a constant, contradicting the inequality above. Therefore, when started from a position $a$ as selected above, the rank process almost surely has no right limit at zero. But, $\lambda(\{a\})>0$. Hence  $(\Phi_t^\nu(X^{\nu,\lambda}(t)))_{t>0}$ cannot be extended to an element in $D([0,\infty))$. 
\end{remark}

\subsection{The Hermite root flow}
\label{sec:hermite entrance}
Here we discuss the special case of the finite free heat flow on the initial polynomial $x^n$. Set

\[
p_{n,t}(x) = T_t\topp n p_{n,0}(x) = \pp{\frac tn}^{n/2}H_n\pp{\sqrt{\frac nt}x},
\]
where $H_n(x)$ is the probabilist's Hermite polynomial of degree $n$. 
If we let $h_{n,1}<\cdots<h_{n,n}$ denote the $n$ roots of $H_n(x)$, then the $n$ roots are now
\[
\vec r_n(t)  = \sqrt{\frac tn}\times \vec h_n, t\ge 0.
\]
The Hermite root flow in this case still satisfies the ODE \eqref{eq:frozen DBm} when restricted to $t\ge \epsilon>0$. 
Moreover, this time
\[
Q\topp n(\vec r_n(t)) = \frac 1t Q\topp n\pp{\frac {\vec h_n}{\sqrt n}},
\] 
and it  has eigenvalues
$-k/(2t)$, $k=0,\dots,n-1$. Consequently, the corresponding $J_{s,t}\topp n$ has eigenvalues $(s/t)^{k/2}, k=0,\dots,n-1$ for $0<s<t$. These facts follow from classical spectral identities for Hermite zeros; see \citet[Section~4.1, equations~(4.5)--(4.7)]{sasaki15perturbations}, where the relevant matrix and its spectral properties are attributed to \citet{ahmed79properties}.

We have defined $X^{\nu,\lambda}$ in Section \ref{sec:classical free Bm}, and let $(X^*(t))_{t\ge 0} \equiv \pp{X^{\delta_0,\delta_0}(t)}_{t\ge0}$  
denote the standard classical free Brownian motion.  The only issue to be dealt with is at $t=0$, where $0$ is a root with multiplicity $n$, and hence we cannot apply the general framework developed earlier directly. For each $n\in\N$, we shall define a Markov process $(X_n^*(t))_{t\ge0}$ through approximation, such that $X_n^*\weakto X^*$ in $D([0,\infty))$. The approximation procedure starts as follows. First, for each $\epsilon>0$, we define $I_n^\epsilon = (I_n^\epsilon(t))_{t\ge\epsilon}$ to be a rank process as before but starting from time $\epsilon>0$ with initial uniform distribution over $\{1,\dots,n\}$ (and the same generator $\calL\topp n_t$ as in \eqref{eq:I generator}). Set 
\[
X_n^\epsilon(t):= r_{n,I_n^\epsilon(t\vee \epsilon)}(t\vee\epsilon), t\ge 0.
\]
The construction of $X^*_n$ is summarized in the following lemma.
\begin{lemma}\label{lem:X_n^epsilon}
The process $(X_n^\epsilon(t))_{t\ge0}$ converges in distribution in $D([0,\infty))$ as $\epsilon\downarrow 0$. If we let $(X_n^*(t))_{t\ge0}$ denote a process with the limiting distribution, then for every $\epsilon>0$ the process $(X_n^*(t\vee\epsilon))_{t\ge0}$ has the same law as $(X_n^\epsilon(t))_{t\ge0}$. 
\end{lemma}
\begin{proof}
We use the fact that, since $I_n^\epsilon(\epsilon)$ is uniformly distributed over $\{1,\dots,n\}$, it has the same uniform law at any time $t>\epsilon$. Therefore, $I_n^\epsilon = (I_n^\epsilon(t))_{t\ge \epsilon}, \epsilon>0$ is a consistent family of processes in the sense that $I_n^\epsilon$ restricted to $[\epsilon',\infty)$ for every $\epsilon'>\epsilon$ has the same law as $I_n^{\epsilon'}$. Therefore, $(X_n^\epsilon)_{\epsilon>0}$ is also consistent. We can thus define a process indexed by $t\in(0,\infty)$, denoted by $(X_n^*(t))_{t\in(0,\infty)}$ with paths in $D((0,\infty))$ satisfying the stated property. Moreover, $\max_{i=1,\dots,n}|r_{n,i}(t)|\le C_n\sqrt t$ for all $t>0$. So we simply set $X_n^*(0) = 0$ by continuity. 
\end{proof}
The next result is the Hermite root flow counterpart for Theorem \ref{thm:2}. Note that we can still couple $X_n^*$ with a rank process $I_n$ restricted to $t\in[\epsilon,\infty)$ for any $\epsilon>0$. We choose to omit this construction.
\begin{proposition}
Let $(X_n^*(t))_{t\ge0}$ denote a process with the limiting law as in Lemma \ref{lem:X_n^epsilon}.
Then, 
\[
(X_n^*(t))_{t\ge0}\weakto (X^*(t))_{t\ge 0}
\]
in $D([0,\infty))$ as $n\to\infty$.
\end{proposition}
\begin{proof}
Fix $\epsilon>0$. Applying Theorem \ref{thm:2} to the tagged-root process indexed by $[\epsilon,\infty)$, and noticing that by construction both $(X_n^{\epsilon}(t))_{t\ge0}$ and $(X^*(t\vee\epsilon))_{t\ge0}$ are constant over $[0,\epsilon]$, we have
\[
(X_n^{\epsilon}(t))_{t\ge 0} \weakto (X^*(t\vee\epsilon))_{t\ge0},
\]
in $D([0,\infty))$. 
It is also straightforward to check that 
\[
(X^*(t\vee\epsilon))_{t\ge0} \weakto (X^*(t))_{t\ge 0}
\]
as $\epsilon\downarrow0$ in $D([0,\infty))$, as the process value at time $t$ is bounded by $2\sqrt t$.  To prove the proposition, it suffices to establish 
\[
\lim_{\epsilon\downarrow0}\limsupn  \proba\pp{\sup_{t\in[0,\epsilon]}\abs{X_n^*(t)-X_n^*(t\vee\epsilon)}>\eta} = 0, \mfa \eta>0.
\]
But, the probability in the expression above is zero for $\epsilon<\eta^2/16$, as $|X_n^*(t)|\le 2\sqrt t$ almost surely. 
\end{proof}

\appendix
\section{Connections to classical free Meixner processes}\label{sec:meixner}
We present a connection between classical free Brownian motion and classical free Meixner processes studied by \citet{bryc05conditional}, where the classical free Meixner processes are introduced as a family of classical Markov processes indexed by parameters $\tau\ge0$ and $\theta\in\R$ with explicit formulas for their transition kernels provided and starting from 0 at $t=0$. 

It is known that all the classical free Meixner processes are classical free Brownian motions in the sense considered here following \citet{anshelevich13generators}. 
Set
 \[
 \omega_{t,\theta} := \mu_{\rmsc}^{\boxplus t}*\delta_\theta, \theta\in\R,
 \] 
 for the translation of measure $\mu_{\rmsc}^{\boxplus t}$ by $\theta$.  \citet{bryc10markov} proved that the classical free Meixner processes are indexed by two parameters, $\tau\ge 0$ and $\theta\in\R$, and in our notation each process has generator $\calA_t^{\omega_{\tau,\theta}}$. We let $M_{\tau,\theta}$ denote such a process. 

What seems to have been missing in the literature is the following connection: when $\tau=0$, $M_{0,\theta}$ has the law of $\theta+X^{\delta_0,\delta_{-\theta}}$, and when $\tau>0$, $M_{\tau,\theta}$ can be viewed as the increment process of a standard classical free Brownian motion at the time $\tau$ starting from $-\theta$. Namely, 
we let
$\spp{\wt X^*_{\tau,-\theta}(t)}_{t\ge 0}$
denote a Markov process with transition kernels  $(\calP^{\delta_0}_{\tau+s,\tau+t})_{0\le s<t<\infty}$ and initial law $\delta_{-\theta}$, and the aforementioned increment process can be represented as
\[
\pp{\theta+\wt X^*_{\tau,-\theta}(t)}_{t\ge0}.
\]
Note that the above representation allows $\tau=0$. We have the following connection. 
\begin{proposition}
For each $\tau\ge 0,\theta\in\R$, 
\equh\label{eq:shift}
\pp{M_{\tau,\theta}(t)}_{t\ge0}\eqd \pp{X^{\omega_{\tau,\theta},\delta_0}(t)}_{t\ge 0} \eqd\pp{\theta+\wt X^*_{\tau,-\theta}(t)}_{t\ge0},
\eque
as random elements in $D([0,\infty))$.
\end{proposition}
\begin{proof}
For $f\in\calW$, recall the definition of $K_f(x,y)$ in \eqref{eq:K_f}. \citet[Corollary 1.7]{bryc14infinitesimal} computed the generator of $M_{\tau,\theta}$ on $\calW$ as
\[
\mathbb A_t^{\tau,\theta}f(x) = \int_\R K_f(x,y)\omega_{\tau+t,\theta}(dy), \mfa f\in \calW.
\]
(For analytic $f$, \citet[Theorem 2.3]{bryc10markov} computed the generator already.)
But, we also have
\[
\mathbb A_t^{\tau,\theta}f(x) = \int_\R K_f(x,y)\pp{\omega_{\tau,\theta}\boxplus \mu_{\rmsc}^{\boxplus t}}(dy) = \calA_t^{\omega_{\tau,\theta}}f(x).
\]
So $M_{\tau,\theta}$ and $X^{\omega_{\tau,\theta},\delta_0}$ have the same generators and the initial law $\delta_0$. Both are Markov processes and they solve the same $(\calA_t^{\omega_{\tau,\theta}},\calW)$-martingale problem. By the uniqueness for the martingale problem, this proves the first equality in law. 

For the second equality in distribution in \eqref{eq:shift}, writing $f_\theta(u):=f(\theta+u)$, we have 
\begin{align*}
K_f(x,\theta+y)&= \frac{f(\theta+y)-f(x)-(\theta+y-x)f'(x)}{(\theta+y-x)^2}\\
& = \frac{f_\theta(y)-f_\theta(x-\theta)- (y-(x-\theta))f_\theta'(x-\theta)}{(y-(x-\theta))^2} = K_{ f_\theta}(x-\theta,y).
\end{align*}
Therefore,  a change of variables gives
\begin{align*}
\calA_t^{\omega_{\tau,\theta}}f(x)
&=\int_\R K_f(x,y)\pp{\delta_\theta* \mu_{\rmsc}^{\boxplus (\tau+t)}}(\d y)=\int_\R K_f(x,\theta+y) \mu_{\rmsc}^{\boxplus (\tau+t)}(\d y) \\
&=\int_\R K_{f_\theta}(x-\theta,y)\mu_{\rmsc}^{\boxplus (\tau+t)}(\d y)
=\calA_{\tau+t}^{\delta_0}f_\theta(x-\theta).
\end{align*}
On the other hand, letting $\wt\calA_t^{\tau,\theta},\wt\calP_{s,t}^{\tau,\theta}$ denote the generator and the transition operator of the process $(\theta+\wt X^*_{\tau,-\theta}(t))_{t\ge0}$, by the definition of $\wt X^*_{\tau,-\theta}$ we have
\begin{align*}
\wt \calP_{t,t+h}^{\tau,\theta}f(x)& = \esp \pp{f\pp{\theta+\wt X_{\tau,-\theta}^*(t+h)}\mmid \theta+\wt X_{\tau,-\theta}^*(t) = x} \\
& = \esp \pp{f_\theta\pp{\wt X_{\tau,-\theta}^*(t+h)}\mmid \wt X_{\tau,-\theta}^*(t) = x-\theta} = \calP^{\delta_0}_{\tau+t,\tau+t+h}f_\theta\pp{x-\theta},
\end{align*}
and hence
\begin{align*}
\wt\calA_t^{\tau,\theta}f(x)& = \lim_{h\downarrow0}\frac{\wt \calP_{t,t+h}^{\tau,\theta}f(x) -f(x)}h\\
& = \lim_{h\downarrow0}\frac{\calP_{\tau+t,\tau+t+h}^{\delta_0}f_\theta(x-\theta) -f_\theta(x-\theta)}h=\calA_{\tau+t}^{\delta_0}f_\theta(x-\theta).
\end{align*}
That is, the Markov process $\theta+\wt X^*_{\tau,-\theta}$ also has the same generator $\calA_t^{\omega_{\tau,\theta}}$. Uniqueness of the solution of the martingale problem proves the second equality in law. 
\end{proof}
Recall that the class consists of five sub-classes of processes as listed below:
\begin{enumerate}[(i)]
\item $\tau=0,\theta=0$: classical free Brownian motion,
\item $\tau=0, \theta\ne 0$: classical free Poisson process,
\item $\tau>0, \theta^2>4\tau$: classical free Pascal process,
\item $\tau>0, \theta^2 = 4\tau$: classical free Gamma process,
\item $\tau>0, \theta^2<4\tau$: classical free Meixner process. 
\end{enumerate}
For example, the classical free Gamma process with parameters $(\tau,\theta)$ has the same law as the increment process of a standard classical free Brownian motion starting at time $\tau$ from a boundary position $-\theta$ (which necessarily takes the value of either $2\sqrt\tau$ or $-2\sqrt\tau$). 

	\bibliographystyle{apalike}
\bibliography{references,references18}

\end{document}